%Paper: q-alg/9412016
%From: chered@math.unc.edu (Ivan Cherednik)
%Date: Fri, 30 Dec 1994 14:58:31 +0500

%AMSTeX
% formatam.tex -- AMSTeX template file
% Version: June 10, 1992
%AMSTEX "MACDONALD'S EVALUATION CONJECTURES AND
%        DIFFERENCE FOURIER TRANSFORM", IVAN CHEREDNIK
\input amstex
%\documentstyle{amams} % input Annals of Mathematics macros.
%%%%%%%%%%%%%%%%%%%%%%%%%%%%%%%%%%%%%%%%%%%%%%%%%%%%%%%%%%%
%% amams.sty: AMSTeX Macros for Articles to be published in
%%
%% Annals of Mathematics
%%
%% Princeton University and the
%% Institute for Advanced Study
%%
%% Published by Princeton University Press
%%%%%%%%%%%%%%%%%%%%%%%%%%%%%%%%%%%%%%%%%%%%%%%%%%%%%%%%%%%

%%%%%%%%%%%%%%%%%%%%%%%%%%%%%%%%%%%%%%%%%%%%%%%%%%%%%%%%%%%
%% Variations on AMSPPT.sty written by Amy Hendrickson
%% TeXnology Inc, Brookline, MA
%% 617 738-8029, amyh@ai.mit.edu
%%%%%%%%%%%%%%%%%%%%%%%%%%%%%%%%%%%%%%%%%%%%%%%%%%%%%%%%%%%

\def\spaces{\space\space\space\space\space\space\space\space\space\space}
\def\spacess{\message{\spaces\spaces\spaces\spaces\spaces\spaces\spaces}}
\spacess
\spacess
\message{Annals of Mathematics Style: Current Version: 1.1. June 10, 1992}
\spacess
\spacess
%%%%%%%%%%%%%%%%%%%%%%%%%%%%%%%%%%%%%%%%%%%%%%%%%%%%%%%%

\catcode`\@=11

\hyphenation{acad-e-my acad-e-mies af-ter-thought anom-aly anom-alies
an-ti-deriv-a-tive an-tin-o-my an-tin-o-mies apoth-e-o-ses apoth-e-o-sis
ap-pen-dix ar-che-typ-al as-sign-a-ble as-sist-ant-ship as-ymp-tot-ic
asyn-chro-nous at-trib-uted at-trib-ut-able bank-rupt bank-rupt-cy
bi-dif-fer-en-tial blue-print busier busiest cat-a-stroph-ic
cat-a-stroph-i-cally con-gress cross-hatched data-base de-fin-i-tive
de-riv-a-tive dis-trib-ute dri-ver dri-vers eco-nom-ics econ-o-mist
elit-ist equi-vari-ant ex-quis-ite ex-tra-or-di-nary flow-chart
for-mi-da-ble forth-right friv-o-lous ge-o-des-ic ge-o-det-ic geo-met-ric
griev-ance griev-ous griev-ous-ly hexa-dec-i-mal ho-lo-no-my ho-mo-thetic
ideals idio-syn-crasy in-fin-ite-ly in-fin-i-tes-i-mal ir-rev-o-ca-ble
key-stroke lam-en-ta-ble light-weight mal-a-prop-ism man-u-script
mar-gin-al meta-bol-ic me-tab-o-lism meta-lan-guage me-trop-o-lis
met-ro-pol-i-tan mi-nut-est mol-e-cule mono-chrome mono-pole mo-nop-oly
mono-spline mo-not-o-nous mul-ti-fac-eted mul-ti-plic-able non-euclid-ean
non-iso-mor-phic non-smooth par-a-digm par-a-bol-ic pa-rab-o-loid
pa-ram-e-trize para-mount pen-ta-gon phe-nom-e-non post-script pre-am-ble
pro-ce-dur-al pro-hib-i-tive pro-hib-i-tive-ly pseu-do-dif-fer-en-tial
pseu-do-fi-nite pseu-do-nym qua-drat-ics quad-ra-ture qua-si-smooth
qua-si-sta-tion-ary qua-si-tri-an-gu-lar quin-tes-sence quin-tes-sen-tial
re-arrange-ment rec-tan-gle ret-ri-bu-tion retro-fit retro-fit-ted
right-eous right-eous-ness ro-bot ro-bot-ics sched-ul-ing se-mes-ter
semi-def-i-nite semi-ho-mo-thet-ic set-up se-vere-ly side-step sov-er-eign
spe-cious spher-oid spher-oid-al star-tling star-tling-ly
sta-tis-tics sto-chas-tic straight-est strange-ness strat-a-gem strong-hold
sum-ma-ble symp-to-matic syn-chro-nous topo-graph-i-cal tra-vers-a-ble
tra-ver-sal tra-ver-sals treach-ery turn-around un-at-tached un-err-ing-ly
white-space wide-spread wing-spread wretch-ed wretch-ed-ly Brown-ian
Eng-lish Euler-ian Feb-ru-ary Gauss-ian Grothen-dieck Hamil-ton-ian
Her-mit-ian Jan-u-ary Japan-ese Kor-te-weg Le-gendre Lip-schitz
Lip-schitz-ian Mar-kov-ian Noe-ther-ian No-vem-ber Rie-mann-ian
Schwarz-schild Sep-tem-ber Za-mo-lod-chi-kov Knizh-nik quan-tum Op-dam
Mac-do-nald Ca-lo-ge-ro Su-ther-land Mo-ser Ol-sha-net-sky  Pe-re-lo-mov }

\Invalid@\nofrills
\Invalid@\usualspace
\newif\ifnofrills@
\def\nofrills@#1#2{\relaxnext@
  \DN@{\ifx\next\nofrills
    \nofrills@true\let#2\relax\DN@\nofrills{\nextii@}%
  \else
    \nofrills@false\def#2{#1}\let\next@\nextii@\fi
\next@}}
\def\usualspace@#1{\ifnofrills@\def\usualspace{#1}\fi}
\def\addto#1#2{\csname \expandafter\eat@\string#1@\endcsname
  \expandafter{\the\csname \expandafter\eat@\string#1@\endcsname#2}}
\newdimen\bigsize@
\def\big@#1#2{{\hbox{$\left#2\vcenter to#1\bigsize@{}%
  \right.\nulldelimiterspace\z@\m@th$}}}
\def\big{\big@\@ne}
\def\Big{\big@{1.5}}
\def\bigg{\big@\tw@}
\def\Bigg{\big@{2.5}}
\def\raggedcenter@{\leftskip\z@ plus.4\hsize \rightskip\leftskip
 \parfillskip\z@ \parindent\z@ \spaceskip.3333em \xspaceskip.5em
 \pretolerance9999\tolerance9999 \exhyphenpenalty\@M
 \hyphenpenalty\@M \let\\\linebreak}
\def\upperspecialchars{\def\ss{SS}\let\i=I\let\j=J\let\ae\AE\let\oe\OE
  \let\o\O\let\aa\AA\let\l\L}
\def\uppercasetext@#1{%
  {\spaceskip1.2\fontdimen2\the\font plus1.2\fontdimen3\the\font
   \upperspecialchars\uctext@#1$\m@th\aftergroup\eat@$}}
\def\uctext@#1$#2${\endash@#1-\endash@$#2$\uctext@}
\def\endash@#1-#2\endash@{%
\uppercase{#1}\if\notempty{#2}--\endash@#2\endash@\fi}
\def\runaway@#1{\DN@{#1}\ifx\envir@\next@
  \Err@{You seem to have a missing or misspelled \string\end#1 ...}%
  \let\envir@\empty\fi}
\newif\iftemp@
\def\notempty#1{TT\fi\def\test@{#1}\ifx\test@\empty\temp@false
  \else\temp@true\fi \iftemp@}

%\comment%%% remove
\font@\tensmc=cmcsc10
\font@\sevenex=cmex7
\font@\sevenit=cmti7
\font@\eightrm=cmr8 % preloaded in plain.tex
\font@\sixrm=cmr6 % preloaded in plain.tex
\font@\eighti=cmmi8     \skewchar\eighti='177 % preloaded
\font@\sixi=cmmi6       \skewchar\sixi='177   % preloaded
\font@\eightsy=cmsy8    \skewchar\eightsy='60 % preloaded
\font@\sixsy=cmsy6      \skewchar\sixsy='60   % preloaded
\font@\eightex=cmex8 %
\font@\eightbf=cmbx8 % preloaded in plain.tex
\font@\sixbf=cmbx6   % preloaded in plain.tex
\font@\eightit=cmti8 % preloaded in plain.tex
\font@\eightsl=cmsl8 % preloaded in plain.tex
\font@\eightsmc=cmcsc10
\font@\eighttt=cmtt8 % preloaded in plain.tex
%\font@\ninerm=cmr9
%\font@\ninei=cmmi9    \skewchar\ninei='177
%\font@\ninesy=cmsy9   \skewchar\ninesy='60
%\font@\nineex=cmex9
%\font@\ninebf=cmbx9
%\font@\nineit=cmti9
%\font@\ninesl=cmsl9
%\font@\ninesmc=cmcsc9
%\font@\ninemsa=msam9
%\font@\ninemsb=msbm9
%\font@\nineeufm=eufm9
%\endcomment%%%

\loadmsam
\loadmsbm
\loadeufm
\UseAMSsymbols

\def\penaltyandskip@#1#2{\relax\ifdim\lastskip<#2\relax\removelastskip
      \ifnum#1=\z@\else\penalty@#1\relax\fi\vskip#2%
  \else\ifnum#1=\z@\else\penalty@#1\relax\fi\fi}
\def\nobreak{\penalty\@M
  \ifvmode\def\penalty@{\let\penalty@\penalty\count@@@}%
  \everypar{\let\penalty@\penalty\everypar{}}\fi}
\let\penalty@\penalty

\def\block{\RIfMIfI@\nondmatherr@\block\fi
       \else\ifvmode\vskip\abovedisplayskip\noindent\fi
        $$\def\endblock{\par\egroup$$}\fi
  \vbox\bgroup\advance\hsize-2\indenti\noindent}
\def\endblock{\par\egroup}

\def\logo@{\baselineskip2pc \hbox to\hsize{\hfil\eightpoint Typeset by
 \AmSTeX}}

%%%%%%%%%%%%%%%%%%%%%%%%%%%%%%%%%%%%%%%%%%%%%%%%%%%%%%%%%%%%%%%
%% Macros for Annals of Mathematics written by Amy Hendrickson
%% TeXnology Inc, Brookline, MA
%% 617 738-8029, amyh@ai.mit.edu
%%%%%%%%%%%%%%%%%%%%%%%%%%%%%%%%%%%%%%%%%%%%%%%%%%%%%%%%%%%%%%%

%% This file includes:
%% 1) Font declarations,
%% 2) Page set up,
%% 3) Title page
%% 4) Section heads,
%% 5) Equation macros, autonumbering equations, etc.,
%% 6) Figure and Table Captions,
%% 7) End matter macros: Bibliography, Appendix, etc.,
%% 8) Footnotes,
%% 9) Theorem type environments
%% 10) Cross-referencing
%% 11) Listing
%% 12) Article and Journal Table of Contents

%%%%%%%%%%%%%%%%%%%%%%%%%%%%%%%%%%%
%% 1) Font declarations,
% Computer Modern fonts

% Small Caps
\font\elevensc=cmcsc10 scaled\magstephalf
\font\tensc=cmcsc10

\font\eightsc=cmcsc10 scaled800

\font\elevenrm=cmr10 scaled \magstephalf%!!!
\font\ninerm=cmr9
\font\eightrm=cmr8
\font\sixrm=cmr6
\font\fiverm=cmr5

\font\eleveni=cmmi10 scaled\magstephalf
\font\ninei=cmmi9
\font\eighti=cmmi8
\font\sixi=cmmi6
\font\fivei=cmmi5
\skewchar\ninei='177 \skewchar\eighti='177 \skewchar\sixi='177
\skewchar\eleveni='177

\font\elevensy=cmsy10 scaled\magstephalf
\font\ninesy=cmsy9
\font\eightsy=cmsy8
\font\sixsy=cmsy6
\font\fivesy=cmsy5
\skewchar\ninesy='60 \skewchar\eightsy='60 \skewchar\sixsy='60
\skewchar\elevensy'60

\font\eighteenbf=cmbx10 scaled\magstep3

\font\twelvebf=cmbx10 scaled \magstep1
\font\elevenbf=cmbx10 scaled \magstephalf
\font\tenbf=cmbx10
\font\ninebf=cmbx9
\font\eightbf=cmbx8
\font\sixbf=cmbx6
\font\fivebf=cmbx5

\font\elevenit=cmti10 scaled\magstephalf
\font\nineit=cmti9
\font\eightit=cmti8

% Fonts for bold math
\font\eighteenmib=cmmib10 scaled \magstep3
\font\twelvemib=cmmib10 scaled \magstep1
\font\elevenmib=cmmib10 scaled\magstephalf
\font\tenmib=cmmib10
\font\eightmib=cmmib10 scaled 800
\font\sixmib=cmmib10 scaled 600

\font\eighteensyb=cmbsy10 scaled \magstep3
\font\twelvesyb=cmbsy10 scaled \magstep1
\font\elevensyb=cmbsy10 scaled \magstephalf
\font\tensyb=cmbsy10
\font\eightsyb=cmbsy10 scaled 800
\font\sixsyb=cmbsy10 scaled 600

\font\elevenex=cmex10 scaled \magstephalf
\font\tenex=cmex10
\font\eighteenex=cmex10 scaled \magstep3

%%%%%%%%%%%%%%%%%%%%%%%%%%%%
%% Font families

\def\elevenpoint{\def\rm{\fam0\elevenrm}%
  \textfont0=\elevenrm \scriptfont0=\eightrm \scriptscriptfont0=\sixrm
  \textfont1=\eleveni \scriptfont1=\eighti \scriptscriptfont1=\sixi
  \textfont2=\elevensy \scriptfont2=\eightsy \scriptscriptfont2=\sixsy
  \textfont3=\elevenex \scriptfont3=\tenex \scriptscriptfont3=\tenex
  \def\bf{\fam\bffam\elevenbf}%
  \def\it{\fam\itfam\elevenit}%
  \textfont\bffam=\elevenbf \scriptfont\bffam=\eightbf
   \scriptscriptfont\bffam=\sixbf
\normalbaselineskip=13.95pt
  \setbox\strutbox=\hbox{\vrule height9.5pt depth4.4pt width0pt\relax}%
  \normalbaselines\rm}

\elevenpoint %%% default fonts and baselineskip

\def\ninepoint{\def\rm{\fam0\ninerm}%
  \textfont0=\ninerm \scriptfont0=\sixrm \scriptscriptfont0=\fiverm
  \textfont1=\ninei \scriptfont1=\sixi \scriptscriptfont1=\fivei
  \textfont2=\ninesy \scriptfont2=\sixsy \scriptscriptfont2=\fivesy
  \textfont3=\tenex \scriptfont3=\tenex \scriptscriptfont3=\tenex
  \def\it{\fam\itfam\nineit}%
  \textfont\itfam=\nineit
  \def\bf{\fam\bffam\ninebf}%
  \textfont\bffam=\ninebf \scriptfont\bffam=\sixbf
   \scriptscriptfont\bffam=\fivebf
\normalbaselineskip=11pt
  \setbox\strutbox=\hbox{\vrule height8pt depth3pt width0pt\relax}%
  \normalbaselines\rm}

\def\eightpoint{\def\rm{\fam0\eightrm}%
  \textfont0=\eightrm \scriptfont0=\sixrm \scriptscriptfont0=\fiverm
  \textfont1=\eighti \scriptfont1=\sixi \scriptscriptfont1=\fivei
  \textfont2=\eightsy \scriptfont2=\sixsy \scriptscriptfont2=\fivesy
  \textfont3=\tenex \scriptfont3=\tenex \scriptscriptfont3=\tenex
  \def\it{\fam\itfam\eightit}%
  \textfont\itfam=\eightit
  \def\bf{\fam\bffam\eightbf}%
  \textfont\bffam=\eightbf \scriptfont\bffam=\sixbf
   \scriptscriptfont\bffam=\fivebf
\normalbaselineskip=12pt
  \setbox\strutbox=\hbox{\vrule height8.5pt depth3.5pt width0pt\relax}%
  \normalbaselines\rm}

%%%%%%%%%%%%%%%%%%%%%%%%%%%%
%% Font families for bold math in title and section heads

\def\eighteenbold{\def\rm{\fam0\eighteenbf}%
  \textfont0=\eighteenbf \scriptfont0=\twelvebf \scriptscriptfont0=\tenbf
  \textfont1=\eighteenmib \scriptfont1=\twelvemib\scriptscriptfont1=\tenmib
  \textfont2=\eighteensyb \scriptfont2=\twelvesyb\scriptscriptfont2=\tensyb
  \textfont3=\eighteenex \scriptfont3=\tenex \scriptscriptfont3=\tenex
  \def\bf{\fam\bffam\eighteenbf}%
  \textfont\bffam=\eighteenbf \scriptfont\bffam=\twelvebf
   \scriptscriptfont\bffam=\tenbf
\normalbaselineskip=20pt
  \setbox\strutbox=\hbox{\vrule height13.5pt depth6.5pt width0pt\relax}%
\everymath {\fam0 }
\everydisplay {\fam0 }
  \normalbaselines\rm}

\def\elevenbold{\def\rm{\fam0\elevenbf}%
  \textfont0=\elevenbf \scriptfont0=\eightbf \scriptscriptfont0=\sixbf
  \textfont1=\elevenmib \scriptfont1=\eightmib \scriptscriptfont1=\sixmib
  \textfont2=\elevensyb \scriptfont2=\eightsyb \scriptscriptfont2=\sixsyb
  \textfont3=\elevenex \scriptfont3=\elevenex \scriptscriptfont3=\elevenex
  \def\bf{\fam\bffam\elevenbf}%
  \textfont\bffam=\elevenbf \scriptfont\bffam=\eightbf
   \scriptscriptfont\bffam=\sixbf
\normalbaselineskip=14pt
  \setbox\strutbox=\hbox{\vrule height10pt depth4pt width0pt\relax}%
\everymath {\fam0 }
\everydisplay {\fam0 }
  \normalbaselines\bf}

%%%%%%%%%%%%%%%%%%%%%%%%%%%%%%%%%%%%%%%%%%%%%%%%%%%%%%%%%
%% 2) Page set up
\hsize=31pc
\vsize=48pc

\parindent=22pt
\parskip=0pt

\widowpenalty=10000
\clubpenalty=10000

\topskip=12pt

\skip\footins=20pt
\dimen\footins=3in % maximum footnote height

\abovedisplayskip=6.95pt plus3.5pt minus 3pt
\belowdisplayskip=\abovedisplayskip

%% Output routine

\voffset=7pt\hoffset= .7in%7pt magstep1

\newif\iftitle%!

\def\amheadline{\iftitle%
\hbox to\hsize{\hss\currannalsline\hss}\else\line{\ifodd\pageno
\hfill\thetitle\hfill\llap{\elevenrm\folio}\else\rlap{\elevenrm\folio}
\hfill\theauthors\hfill\fi}\fi}

\headline={\amheadline}%!!!
\footline={\global\titlefalse}
%\output={\bindingoffset\plainoutput}

%%%%%%%%%%%%%%%%%%%%%%%%%%%%%%%%%%%%%%%
% 3) Title page

 %#1= Volume number, #2=year of publication
\def\annalsline#1#2{\vfill\eject
\ifodd\pageno\else % first page of article on right.
\line{\hfill}
\vfill\eject\fi
\global\titletrue
\def\currannalsline{\eightrm %Annals of Mathematics,%ANNALS
{\eightbf#1} (#2), \thepages}}

\def\titleheadline#1{\def\one{#1}\ifx\one\empty\else
\def\thetitle{{%\frenchspacing%
\let\\ \relax\eightsc\uppercase{#1}}}\fi}

\newif\ifshort

\let\shorttitle\titleheadline

\def\onpages#1#2{\def\thepages{#1--#2}}

\def\thismuchskip[#1]{\vskip#1pt}
\def\ilook{\ifx\next[ \let\go\thismuchskip\else
\let\go\relax\vskip1pt\fi\go}

\def\institution#1{\def\theinstitutions{\vbox{\baselineskip10pt
\def\and{{\eightrm and }}
\def\\{\futurelet\next\ilook}\eightsc #1}}}
\let\institutions\institution

\newwrite\auxfile

\def\startingpage#1{\def\one{#1}\ifx\one\empty\global\pageno=1\else
\global\pageno=#1\fi
\theoremcount=0 \eqcount=0 \sectioncount=0
\openin1 \jobname.aux \ifeof1
\onpages{#1}{???}
\else\closein1 \relax\input \jobname.aux
\onpages{#1}{\lastpage}
\fi\immediate\openout\auxfile=\jobname.aux
}

\def\endarticle{\ifRefsUsed\global\RefsUsedfalse%
\else\vskip21pt\theinstitutions%
\nobreak\vskip8pt
%\vbox{\thereceived\therevised}%
\fi%
\write\auxfile{\string\def\string\lastpage{\the\pageno}}}

\outer\def\bye{\endarticle\par \vfill \supereject \end}

% variation on code from amsspt.sty ==>
\def\document{\let\fontlist@\relax\let\alloclist@\relax
 \elevenpoint}%%% add for annals!!!

% <=== end of code varied from amsppt.sty

\newif\ifacks
\long\def\acknowledgements#1{\def\one{#1}\ifx\one\empty\else
\vskip-\baselineskip%
\global\ackstrue\footnote{\ \unskip}{*#1}\fi}

\def\title#1{\titleheadline{#1}
\vbox to80pt{\vfill
\baselineskip=18pt
\parindent=0pt
\overfullrule=0pt
\hyphenpenalty=10000
\everypar={\hskip\parfillskip\relax}
\hbadness=10000
\def\\ {\vskip1sp}
\eighteenbold#1\vskip1sp}}

\newif\ifauthor

\def\author#1{\vskip11pt
\hbox to\hsize{\hss\tenrm By \tensc#1\ifacks\global\acksfalse*\fi\hss}
\ifshort\else\xdef\theauthors{{\eightsc\uppercase{#1}}}\fi%
\vskip21pt\global\authortrue\everypar={\global\authorfalse\everypar={}}}

\def\twoauthors#1#2{\vskip11pt
\hbox to\hsize{\hss%
\tenrm By \tensc#1 {\tenrm and} #2\ifacks\global\acksfalse*\fi\hss}
\ifshort\else\xdef\theauthors{{\eightsc\uppercase{#1 and #2}}}\fi%
\vskip21pt
\global\authortrue\everypar={\global\authorfalse\everypar={}}}

%%%%%%%%%%%%%%%%%%%%%%%%%%%%%%%%
%% 4) Section heads, counters

\newcount\theoremcount
\newcount\sectioncount
\newcount\eqcount

\newif\ifspecialnumon

\def\eqnumber=#1 {\global\eqcount=#1 \global\advance\eqcount by-1\relax}
\def\sectionnumber=#1 {\global\sectioncount=#1
\global\advance\sectioncount by-1\relax}
\def\proclaimnumber=#1 {\global\theoremcount=#1
\global\advance\theoremcount by-1\relax}

\newif\ifsection
\newif\ifsubsection

\def\intro{\global\authorfalse%                           !!! remove
\centerline{\bf Introduction}\everypar={}\vskip6pt}

\def\elevenboldmath#1{$#1$\egroup}
\def\mathbold{\hbox\bgroup\elevenbold\elevenboldmath}

\def\section#1{\global\theoremcount=0
\global\eqcount=0
\ifauthor\global\authorfalse\else%
\vskip18pt plus 18pt minus 6pt\fi%
{\parindent=0pt
\everypar={\hskip\parfillskip}%            !!! remove
\def\\ {\vskip1sp}\elevenpoint\bf%
\ifspecialnumon\global\specialnumonfalse$\rm\spnum$%
\gdef\sectnum{$\rm\spnum$}%
\else\interlinepenalty=10000%
\global\advance\sectioncount by1\relax\the\sectioncount%
\gdef\sectnum{\the\sectioncount}%
\fi. \hskip6pt#1%                          !!!add }} and stop here
\vrule width0pt depth12pt}
\hskip\parfillskip%\break%!
\global\sectiontrue%
\everypar={\global\sectionfalse\global\interlinepenalty=0\everypar={}}%
\ignorespaces

}

%%%%%%%%%%%%%%%%%%%%%%%%%%%%%%%%
%% 5) Equation Macros

\newif\ifspequation

\let\eqno\leqno %automatic left side equation numbers %%!!!remove l-eqno

\newif\ifineqalignno
\let\saveleqalignno\leqalignno                        %%!!!remove l-eqno
\def\leqalignno{\let\eqnu\Eeqnu\saveleqalignno}

\let\eqalignno\leqalignno

\def\sectandeqnum{%
\ifspecialnumon\global\specialnumonfalse
$\rm\spnum$\gdef\eqnum{{$\rm\spnum$}}\else\global\firstlettertrue
\global\advance\eqcount by1
\ifappend\applett\else\the\sectioncount\fi.%
\the\eqcount
\xdef\eqnum{\ifappend\applett\else\the\sectioncount\fi.\the\eqcount}\fi}

\def\eqnu{\leqno{\hbox{\elevenrm\ifspequation\else(\fi\sectandeqnum
\ifspequation\global\spequationfalse\else)\fi}}}      %!!! l-eqno

\def\Speqnu{\global\setbox\leqnobox=\hbox{\elevenrm
\ifspequation\else%
(\fi\sectandeqnum\ifspequation\global\spequationfalse\else)\fi}}

\def\Eeqnu{\hbox{\elevenrm
\ifspequation\else%
(\fi\sectandeqnum\ifspequation\global\spequationfalse\else)\fi}}

\newif\iffirstletter
\global\firstlettertrue
\def\eqletter#1{\global\specialnumontrue\iffirstletter\global\firstletterfalse
\global\advance\eqcount by1\fi
\gdef\spnum{\the\sectioncount.\the\eqcount#1}\eqnu}

%%% Split math
\newbox\leqnobox
\def\outsideeqnu#1{\global\setbox\leqnobox=\hbox{#1}}

\def\eatone#1{}

%% Vertically centers equation number.
\def\dosplit#1#2{\vskip-.5\abovedisplayskip
\setbox0=\hbox{$\let\eqno\outsideeqnu%
\let\eqnu\Speqnu\let\leqno\outsideeqnu#2$}%
\setbox1\vbox{\noindent\hskip\wd\leqnobox\ifdim\wd\leqnobox>0pt\hskip1em\fi%
$\displaystyle#1\mathstrut$\hskip0pt plus1fill\relax
\vskip1pt
\line{\hfill$\let\eqnu\eatone\let\leqno\eatone%
\displaystyle#2\mathstrut$\ifmathqed~~\qed\fi}}%
\copy1
\ifvoid\leqnobox
\else\dimen0=\ht1 \advance\dimen0 by\dp1
\vskip-\dimen0
\vbox to\dimen0{\vfill
\hbox{\unhbox\leqnobox}
\vfill}
\fi}

\everydisplay{\lookforbreak}

\long\def\lookforbreak #1$${\def\mathone{#1}
\expandafter\testforbreak\mathone\splitmath @}

\def\testforbreak#1\splitmath #2@{\def\mathtwo{#2}\ifx\mathtwo\empty%
#1$$%
\ifmathqed\vskip-\belowdisplayskip
\setbox0=\vbox{\let\eqno\relax\let\eqnu\relax$\displaystyle#1$}%
\vskip-\ht0\vskip-3.5pt\hbox to\hsize{\hfill\qed}
\vskip\ht0\vskip3.5pt\fi
\else$$\vskip-\belowdisplayskip
\vbox{\dosplit{#1}{\let\eqno\eatone
\let\splitmath\relax#2}}%
\nobreak\vskip.5\belowdisplayskip
\noindent\ignorespaces\fi}

%% Proof box to be used when proof ends with equation.

\newif\ifmathqed

%%%%%%%%%%%%%%%%%%%%%%%%%%%%%
%% \mtable, Math table to make binary table easily

%% Use:
% \mtable
% &n_1&n_2&n_3&n_4&n_5&n_6\cr
% \Delta_1&M_3&M_2&0&0&0&0\cr
% \Delta_2&0&0&M_1&M_3&0&0\cr
% \endmtable

\newcount\linenum
\newcount\colnum

%++
\def\spline{\omit&\multispan{\the\colnum}{\hrulefill}\cr}
\def\colcounter{\ifnum\linenum=1\global\advance\colnum by1\fi}

\def\everyline{\noalign{\global\advance\linenum by1\relax}%
\ifnum\linenum=2\spline\fi}

\def\mtable{\bgroup\offinterlineskip
\everycr={\everyline}\global\linenum=0
\halign\bgroup\vrule height 10pt depth 4pt width0pt
\hfill$##$\hfill\hskip6pt\ifnum\linenum>1
\vrule\fi&&\colcounter\hskip12pt\hfill$##$\hfill\hskip12pt\cr}

\def\endmtable{\crcr\egroup\egroup}

%%%%%%%%%%%%%%%%%%%%%%%%%%%%%
% Array

%% Will work in math or in text, will be in math mode inside array.
%% For each column desired supply
%% r, l, or c, for right, left, or center orientation of that column.
%% End each line with \\.

%% To use:
%  \array ccc*
%  x_s\leq a_1\\
%  a_s<x_s^s<b_s\\
%  x_s\geq a_1
%  \endarray

\def\xast{*}
\newcount\intable
\newcount\mathcol
\newcount\savemathcol
\newcount\topmathcol
\newdimen\arrayhspace
\newdimen\arrayvspace

\arrayhspace=8pt % horizontal space between columns, (half this width
                 %  will horizontally precede and follow the array)
\arrayvspace=12pt % vertical space between lines

\newif\ifdollaron

\def\mathalign#1{\def\arg{#1}\ifx\arg\xast%
\let\go\relax\else\let\go\mathalign%
\global\advance\mathcol by1 %
\global\advance\topmathcol by1 %
\expandafter\def\csname  mathcol\the\mathcol\endcsname{#1}%
\fi\go}

\def\arraypickapart#1]#2*{\if#1c \ifmmode\vcenter\else
\global\dollarontrue$\vcenter\fi\else%
\if#1t\vtop\else\if#1b\vbox\fi\fi\fi\bgroup%
\def\one{#2}}

\def\arraystrut{\vrule height .7\arrayvspace depth .3\arrayvspace width 0pt}

\def\array#1#2*{\def\firstarg{#1}%
\if\firstarg[ \def\two{#2} \expandafter\arraypickapart\two*\else%
\ifmmode\vcenter\else\vbox\fi\bgroup \def\one{#1#2}\fi%
\global\everycr={\noalign{\global\mathcol=\savemathcol\relax}}%
\def\\ {\cr}%
\global\advance\intable by1 %
\ifnum\intable=1 \global\mathcol=0 \savemathcol=0 %
\else \global\advance\mathcol by1 \savemathcol=\mathcol\fi%
\expandafter\mathalign\one*%
\mathcol=\savemathcol %
\halign\bgroup&\hskip.5\arrayhspace\arraystrut%
\global\advance\mathcol by1 \relax%
\expandafter\if\csname mathcol\the\mathcol\endcsname r\hfill\else%
\expandafter\if\csname mathcol\the\mathcol\endcsname c\hfill\fi\fi%
$\displaystyle##$%
\expandafter\if\csname mathcol\the\mathcol\endcsname r\else\hfill\fi\relax%
\hskip.5\arrayhspace\cr}

\def\endarray{\crcr\egroup\egroup%
\global\mathcol=\savemathcol %
\global\advance\intable by -1\relax%
\ifnum\intable=0 %
\ifdollaron\global\dollaronfalse $\fi
\loop\ifnum\topmathcol>0 %
\expandafter\def\csname  mathcol\the\topmathcol\endcsname{}%
\global\advance\topmathcol by-1 \repeat%
\global\everycr={}\fi%
}

\def\big#1{{\hbox{$\left#1\vbox to 10pt{}\right.\n@space$}}}
\def\Big#1{{\hbox{$\left#1\vbox to 13pt{}\right.\n@space$}}}
\def\bigg#1{{\hbox{$\left#1\vbox to 16pt{}\right.\n@space$}}}
\def\Bigg#1{{\hbox{$\left#1\vbox to 19pt{}\right.\n@space$}}}

%%%%%%%%%%%%%%%%%%%%%%%%%%%%%%%%%%%%%%%%%%%%%%%%%%%%%%%%%%%%%%%%
% 6) Figure and Table Captions.

\def\figcaption#1#2#3{\topinsert
\vskip4pt %<===topadjust to match height of ascenders on opposing page.
\vbox to#3{\vfill}\vskip1sp
\setbox0=\hbox{\eightsc Figure #1.\hskip12pt\eightpoint #2}
\ifdim\wd0>\hsize
\noindent\eightsc Figure #1.\hskip12pt\eightpoint #2
\else
\centerline{\eightsc Figure #1.\hskip12pt\eightpoint #2}
\fi
\vskip16pt
\endinsert}

\def\wfig#1#2#3{\topinsert
\vskip4pt %<===topadjust to match height of ascenders on opposing page.
\hbox to\hsize{\hss\vbox{\hrule height .25pt width #3
\hbox to #3{\vrule width .25pt height #2\hfill\vrule width .25pt height #2}
\hrule height.25pt}\hss}
\vskip1sp
\centerline{\eightsc Figure #1}
\vskip16pt
\endinsert}

\def\wfigcaption#1#2#3#4{\topinsert
\vskip4pt %<===topadjust to match height of ascenders on opposing page.
\hbox to\hsize{\hss\vbox{\hrule height .25pt width #4
\hbox to #4{\vrule width .25pt height #3\hfill\vrule width .25pt height #3}
\hrule height.25pt}\hss}
\vskip1sp
\setbox0=\hbox{\eightsc Figure #1.\hskip12pt\eightpoint\rm #2}
\ifdim\wd0>\hsize
\noindent\eightsc Figure #1.\hskip12pt\eightpoint\rm #2\else
\centerline{\eightsc Figure #1.\hskip12pt\eightpoint\rm #2}\fi
\vskip16pt
\endinsert}

\def\tabcaption#1#2{\vskip6pt
\setbox0=\hbox{\eightsc Table #1.\hskip12pt\eightpoint #2}
\ifdim\wd0>\hsize
\noindent\eightsc Table #1.\hskip12pt\eightpoint #2
\else
\centerline{\eightsc Table #1.\hskip12pt\eightpoint #2}
\fi
\vskip6pt}

\def\endinsert{\egroup\if@mid\dimen@\ht\z@\advance\dimen@\dp\z@
\advance\dimen@ 12\p@\advance\dimen@\pagetotal\ifdim\dimen@ >\pagegoal
\@midfalse\p@gefalse\fi\fi\if@mid\smallskip\box\z@\bigbreak\else
\insert\topins{\penalty 100 \splittopskip\z@skip\splitmaxdepth\maxdimen
\floatingpenalty\z@\ifp@ge\dimen@\dp\z@\vbox to\vsize {\unvbox \z@
\kern -\dimen@ }\else\box\z@\nobreak\smallskip\fi}\fi\endgroup}

\def\pagecontents{
\ifvoid\topins \else\iftitle\else
\unvbox \topins \fi\fi \dimen@ =\dp \@cclv \unvbox
\@cclv
\ifvoid\topins\else\iftitle\unvbox\topins\fi\fi
\ifvoid \footins \else \vskip \skip \footins \footnoterule
\unvbox \footins \fi \ifr@ggedbottom \kern -\dimen@ \vfil \fi}

%%%%%%%%%%%%%%%%%%%%%%%%%%%%%%%%%%%%%%%%%%%%%%%%%%%%%%%%%%%%%%%%
% 7) End Matter

\newif\ifappend

\def\appendix#1#2{\def\applett{#1}\def\two{#2}%
\global\appendtrue
\global\theoremcount=0
\global\eqcount=0
\vskip18pt plus 18pt
\vbox{\parindent=0pt
\everypar={\hskip\parfillskip}
\def\\ {\vskip1sp}\elevenbold Appendix%
\ifx\applett\empty\gdef\applett{A}\ifx\two\empty\else.\fi%
\else\ #1.\fi\hskip6pt#2\vskip12pt}%
\global\sectiontrue%
\everypar={\global\sectionfalse\everypar={}}\nobreak\ignorespaces}

\newif\ifRefsUsed
\long\def\references{\global\RefsUsedtrue\vskip21pt
\theinstitutions
\global\everypar={}\global\bibnum=0
\vskip20pt\goodbreak\bgroup
\vbox{\centerline{\eightsc References}\vskip6pt}%
\ifdim\maxbibwidth>0pt
\leftskip=\maxbibwidth%
\parindent=-\maxbibwidth%
\else
\leftskip=18pt%
\parindent=-18pt%
\fi
\ninepoint
\frenchspacing
\nobreak\ignorespaces\everypar={\amref}%
}

\def\endreferences{\vskip1sp\egroup\global\everypar={}%
\nobreak\vskip8pt\vbox{\thereceived\therevised}
}

\newcount\bibnum

\def\amref#1 {\global\advance\bibnum by1%
\immediate\write\auxfile{\string\expandafter\string\def\string\csname
\space #1croref\string\endcsname{[\the\bibnum]}}%
\leavevmode\hbox to18pt{\hbox to13.2pt{\hss[\the\bibnum]}\hfill}}

\def\bibline{\hbox to30pt{\hrulefill}\/\/}

\def\name#1{{\eightsc#1}}

\newdimen\maxbibwidth
\def\AuthorRefNames [#1] {%
\immediate\write\auxfile{\string\def\string\cite\string##1{[\string##1]}}

\def\amref{\spamref}
\setbox0=\hbox{[#1] }\global\maxbibwidth=\wd0\relax}

\def\spamref[#1] {\leavevmode\hbox to\maxbibwidth{\hss[#1]\hfill}}

%%%%%%%%%%%%%%%%%%%%%%%%%%%%%%%%%%%%%%%%%%%%%%%%%%%%%%%%%%%%%%%%
%% 8) Footnotes

\def\footnoterule{\kern-3pt\hrule width1in height.5pt\kern2.5pt}

\def\footnote#1#2{%
\plainfootnote{#1}{{\eightpoint\normalbaselineskip11pt
\normalbaselines#2}}}

\def\vfootnote#1{%
\insert \footins \bgroup \eightpoint\baselineskip11pt
\interlinepenalty \interfootnotelinepenalty
\splittopskip \ht \strutbox \splitmaxdepth \dp \strutbox \floatingpenalty
\@MM \leftskip \z@skip \rightskip \z@skip \spaceskip \z@skip
\xspaceskip \z@skip
{#1}$\,$\footstrut \futurelet \next \fo@t}

%%%%%%%%%%%%%%%%%%%%%%%%%%%%%%%%%%%%%%%%%%%%%%%%%%%%%%%%%%%%%%%%
%% 9) Theorem type environments

\newif\iffirstadded
\newif\ifadded

\def\addedlett{}

\def\alltheoremnums{%
\ifspecialnumon\global\specialnumonfalse
\ifadded\global\addedfalse
\iffirstadded\global\firstaddedfalse
\global\advance\theoremcount by1 \fi
\ifappend\applett\else\the\sectioncount\fi.\the\theoremcount\addedlett%
\xdef\theoremnum{\ifappend\applett\else\the\sectioncount\fi.%
\the\theoremcount\addedlett}%
\else$\rm\spnum$\def\theoremnum{{$\rm\spnum$}}\fi%
\else\global\firstaddedtrue
\global\advance\theoremcount by1
\ifappend\applett\else\the\sectioncount\fi.\the\theoremcount%
\xdef\theoremnum{\ifappend\applett\else\the\sectioncount\fi.%
\the\theoremcount}\fi}

\def\allcorolnums{%
\ifspecialnumon\global\specialnumonfalse
\ifadded\global\addedfalse
\iffirstadded\global\firstaddedfalse
\global\advance\corolcount by1 \fi
\the\corolcount\addedlett%
\else$\rm\spnum$\def\corolnum{$\rm\spnum$}\fi%
\else\global\advance\corolcount by1
\the\corolcount\fi}

%% use for Theorem, Corollary, Lemma, Proposition, Demonstration and similar.

\newcount\corolcount
\def\xcorol{Corollary}
\def\xtheorem{Theorem}
\def\xmaintheorem{Main Theorem}

\newif\ifthtitle

\let\saverparen)
\let\savelparen(
\def\rmparenl{{\rm(}}
\def\rmparenr{{\rm\/)}}
{
\catcode`(=13
\catcode`)=13
\gdef\makeparensRM{\catcode`(=13\catcode`)=13\let(=\rmparenl%
\let)=\rmparenr%
\everymath{\let(\savelparen%
\let)\saverparen}%
\everydisplay{\let(\savelparen%
\let)\saverparen\lookforbreak}}}

\medskipamount=8pt plus.1\baselineskip minus.05\baselineskip

\def\rmtext#1{\hbox{\rm#1}}

\def\proclaim#1{\vskip-\lastskip
\def\one{#1}\ifx\one\xtheorem\global\corolcount=0\fi
\ifsection\global\sectionfalse\vskip-6pt\fi
\medskip
{\elevensc#1}%
\ifx\one\xmaintheorem\global\corolcount=0
\gdef\theoremnum{Main Theorem}\else%
\ifx\one\xcorol\ \allcorolnums\else\ \alltheoremnums\fi\fi%
\ifthtitle\ \global\thtitlefalse{\rm(\thethtitle)}\fi.%
\hskip1em\bgroup\let\text\rmtext\makeparensRM\it\ignorespaces}

\def\nonumproclaim#1{\vskip-\lastskip
\def\one{#1}\ifx\one\xtheorem\global\corolcount=0\fi
\ifsection\global\sectionfalse\vskip-6pt\fi
\medskip
{\elevensc#1}.\ifx\one\xmaintheorem\global\corolcount=0
\gdef\theoremnum{Main Theorem}\fi\hskip.5pc%
\bgroup\it\makeparensRM\ignorespaces}

\def\endproclaim{\egroup\medskip}

%% Use demo for Proof, Proof of, Definition, Example,
%% Remark, Case, Subcase, Conjecture, Note, Notation,
%% Convention, Construction and Step.
%% Any other use for demo will format similar to `Proof.'

\def\xproof{Proof}
\def\xremark{Remark}
\def\xcase{Case}
\def\xsubcase{Subcase}
\def\xconjecture{Conjecture}
\def\xstep{Step}
\def\xof{of}

\def\deconstruct#1 #2 #3 #4 #5 @{\def\one{#1}\def\two{#2}\def\three{#3}%
\def\four{#4}%
\ifx\two\empty #1\else%
\ifx\one\xproof%
\ifx\two\xof%
  \ifx\three\xcorol Proof of Corollary \rm#4\else%
     \ifx\three\xtheorem Proof of Theorem \rm#4\else\xone\fi%
  \fi\fi%
\else\xone\fi\fi.}

\def\pickup#1 {\def\this{#1}%
\ifx\this\xproof\global\let\go\demoproof
\global\let\enddemo\endproof\else
\ifx\this\xremark\global\let\go\demoremark\else
\ifx\this\xcase\global\let\go\demostep\else
\ifx\this\xsubcase\global\let\go\demostep\else
\ifx\this\xconjecture\global\let\go\demostep\else
\ifx\this\xstep\global\let\go\demostep\else
\global\let\go\demoproof\fi\fi\fi\fi\fi\fi}

\newif\ifnonum
\def\demo#1{\vskip-\lastskip
\ifsection\global\sectionfalse\vskip-6pt\fi
\def\one{#1 }\def\two{#1*}%
\setbox0=\hbox{\expandafter\pickup\one}\expandafter\go\two}

\def\numbereddemo#1{\vskip-\lastskip
\ifsection\global\sectionfalse\vskip-6pt\fi
\def\two{#1*}%
\expandafter\demoremark\two}

\def\demoproof#1*{\medskip\def\xone{#1}
{\ignorespaces\it\expandafter\deconstruct\xone {} {} {} {} {} @%
\unskip\hskip6pt}\rm\ignorespaces}

\def\demoremark#1*{\medskip
{\it\ignorespaces#1\/} \ifnonum\global\nonumtrue\else
 \alltheoremnums\unskip.\fi\hskip1pc\rm\ignorespaces}

\def\demostep#1 #2*{\vskip4pt
{\it\ignorespaces#1\/} #2.\hskip1pc\rm\ignorespaces}

\def\enddemo{\medskip}

\def\endproof{\ifmathqed\global\mathqedfalse\medskip\else
\parfillskip=0pt~~\hfill\qed\medskip
\fi\global\parfillskip0pt plus 1fil\relax
\gdef\enddemo{\medskip}}

\def\qed{\vbox{\hrule\hbox{\vrule height6pt\hskip6pt\vrule}\hrule}}

%% Proof box to be used in a \proclaim{}...\endproclaim environment

\def\proofbox{\parfillskip=0pt~~\hfill\qed\vskip1sp\parfillskip=
0pt plus 1fil\relax}

%%%%

%%%%%%%%%%%%%%%%%%%%%
%% 10) CrossRefs

%%% Generic crossreferencing
%%% to use: \label\nameoflabel* (will give the page number when referenced)

% Commands to access current state of counter, for cross-referencing
% \sectnum
% \theoremnum
% \eqnum

%%% You can make another definition that includes counters and/or the
%%% page number and access this information as the second argument:
%%% \label\yourlabelname[2.13]*

%%% Since this method of cross-referencing relies
%%% on an auxiliary file, the first time you tex the file
%%% you will get `??' when you write \ref\nameoflabel.
%%% When you TeX the file the second time the auxiliary file
%%% will be input and \ref\nameoflabel will produce the cross-ref.

\def\stripbs#1#2*{\def\one{#2}}

\def\emptyspace{ }
\def\nextthing{}
\def\newline{***}
\def\eatone#1{ }

\def\lookatspace#1{\ifcat\noexpand#1\ \else%
\gdef\nextthing{}\xdef\next{#1}%
\ifx\next\emptyspace%
\let\nextthing\emptyspace\else\ifx\next\newline%
\gdef\nextthing{\eatone}\fi\fi\fi\egroup\nextthing#1}

{\catcode`\^^M=\active%
\gdef\spacer{\bgroup\catcode`\^^M=\active%
\let^^M=\newline\obeyspaces\lookatspace}}

\def\ref#1{\seeifdefined{#1}\expandafter\csname\one\endcsname\spacer}

\def\cite#1{\expandafter\ifx\csname#1croref\endcsname\relax[??]\else
\csname#1croref\endcsname\fi\spacer}

%% for testing in \label and \ref to see if term already labeled.

\def\seeifdefined#1{\expandafter\stripbs\string#1croref*%
\crorefdefining{#1}}

\newif\ifcromessage
\global\cromessagetrue

\def\crorefdefining#1{\ifdefined{\one}{}
{\ifcromessage\global\cromessagefalse%
\message{\spaces\spaces\spaces\spaces\spaces\spaces\spaces}%
\message{<Undefined reference.}%
\message{Please TeX file once more to have accurate cross-references.>}%
\message{\spaces\spaces\spaces\spaces\spaces\spaces\spaces}\fi[??]}}

\def\label#1#2*{\gdef\ctest{#2}%
\xdef\currlabel{\string#1croref}
\expandafter\seeifdefined{#1}%
\ifx\empty\ctest%
\xdef\labelnow{\write\auxfile{\noexpand\def\currlabel{\the\pageno}}}%
\else\xdef\labelnow{\write\auxfile{\noexpand\def\currlabel{#2}}}\fi%
\labelnow}

\def\ifdefined#1#2#3{\expandafter\ifx\csname#1\endcsname\relax%
#3\else#2\fi}

%%%%%%%%%%%%%%%%%%%%%
%% 11) Listing

%% To use with asterisks:

%%%%%%%%%%%%%%%%%%%%%%
%% 12) Article and Journal Table of Contents

\def\articlecontents{
\vskip20pt\centerline{\bf Table of Contents}\everypar={}\vskip6pt
\bgroup \leftskip=3pc \parindent=-2pc
\def\item##1{\vskip1sp\indent\hbox to2pc{##1.\hfill}}}

\def\endcontents{\vskip1sp\leftskip=0pt\egroup}

\def\journalcontents{\vfill\eject
\def\currannalsline{\hfill}
\global\titletrue
\vglue3.5pc
\centerline{\tensc\hskip12pt TABLE OF CONTENTS}\everypar={}\vskip30pt
\bgroup \leftskip=34pt \rightskip=-12pt \parindent=-22pt
  \def\\ {\vskip1sp\noindent}
\def\pagenum##1{\unskip\parfillskip=0pt\dotfill##1\vskip1sp
\parfillskip=0pt plus 1fil\relax}
\def\name##1{{\tensc##1}}}

%% default values

\institution{}
\onpages{0}{0}
\def\lastpage{???}
\def\thetitle{Title ???}
\def\theauthors{Authors ???}
\def\thereceived{}
\def\therevised{}

\gdef\split{\relaxnext@\ifinany@\let\next\insplit@\else
 \ifmmode\ifinner\def\next{\onlydmatherr@\split}\else
 \let\next\outsplit@\fi\else
 \def\next{\onlydmatherr@\split}\fi\fi\let\eqnu\xspliteqnu\next}

\gdef\align{\relaxnext@\ifingather@\let\next\galign@\else
 \ifmmode\ifinner\def\next{\onlydmatherr@\align}\else
 \let\next\align@\fi\else
 \def\next{\onlydmatherr@\align}\fi\fi\let\eqnu\xspliteqnu\next}

\def\spliteqnu{{\tenrm\sectandeqnum}\relax}

\def\xspliteqnu{\tag\spliteqnu}

\catcode`@=12

\document

%-------------- Publisher's entries --------------------
\annalsline{December}{1994}
%\line{\hfil Devoted to Dani\"el Anton}
\startingpage{1}     %numeration
%%\received{??}
%%\revised{??}

\comment
\nopagenumbers
\headline{\ifnum\pageno=1\hfil\else \rightheadline\fi}
%{\ifodd\pageno\rightheadline \else \leftheadline\fi}\fi}%%after else
\def\rightheadline{\hfil\eightit
%! Running title (odd page)
The Macdonald conjecture
\quad\eightrm\folio}

\voffset=2\baselineskip
\endcomment

%\magnification=\magstep1

%--------------- Author macros ---------------
%                   MACROS
%
%                                 AUX
%
%
%
%                      endaux
%

\def\for{\  \hbox{ for } \ }
\def\if{ \ \hbox{ if } \ }
\def\when{ \ \hbox{ when } \ }
\def\where{\  \hbox{ where } \ }
\def\and{\  \hbox{ and } \ }

\def\equal{\buildrel  def \over =}

\def\la{\lambda}

\def\om{\omega}

\def\th{\theta}
\def\al{\alpha}

\def\ga{\gamma}
\def\ep{\epsilon}

\def\de{\delta}

\def\Ga{\Gamma}
\def\ze{\zeta}

    %from copy, ell

\def\vph{\varphi}

\def\tal{\tilde{\alpha}}

\def\tw{\tilde w}

\def\tz{\tilde z}
\def\tb{\tilde b}

\def\hH{\hat{H}}

\def\hY{\hat{Y}}

\def\hT{\hat{T}}

\def\hw{\hat{w}}

\def\hv{\hat{v}}

\def\C{\bold{C}}
\def\R{\bold{R}}

\def\Z{\bold{Z}}

\def\one{\bold{1}}

\def\0{\bold{0}}

\def\C{\hbox{\bf C}}

%\def\bs{\hbox{\bf S}}          %ell
% macdonald

\def\l{\Cal{L}}
\def\m{\Cal{M}}

\def\p{\Cal{P}}

\def\y{\Cal{Y}}

\def\x{\Cal{X}}
\def\g{\Cal{G}}

\def\w{\Cal{W}}

\font\germ=eufb10 at 12pt
%\font\germm=eufb9 at 12pt
%\font\germ=eufm9 at 12pt
\def\goth#1{\hbox{\germ #1}}

\def\TT{\goth{T}}
\def\HH{\goth{H}}

\font\smm=msbm10 at 12pt
\def\symbol#1{\hbox{\smm #1}}
\def\lsmash{{\symbol n}}

%endmacros

%------------------------------------------------------------------
%-------------- Author entries --------------------

%\comment                               %to remove the title
\title
{Macdonald's Evaluation Conjectures\\
and Difference Fourier Transform
}
 %Article title
\shorttitle{ Evaluation Conjectures}
 % Shortened version for headline title

% Acknowledgements: Please enter all acknowledgements here.
\acknowledgements{
Partially supported by NSF grant DMS--9301114}

% Please uncomment and use appropriate command:
\author{ Ivan Cherednik}
%\twoauthors{}{}
%\authors{}% Separate each author with a comma and a space.

% Institution:
\institutions{
University of North Carolina at Chapel Hill,
Chapel Hill, N.C. 27599-3250
\\ Internet: chered\@math.unc.edu
}

%\endcomment                              %to remove the title
%-------------- Article Text--------------------

\intro %(Optional, Introduction)
%
%
%
%                        INTRO
%
%
%{\bf 0. Introduction.}
\vfil
 Generalizing
the characters of compact simple Lie groups,
Ian Macdonald introduced in [M1,M2] and other works
remarkable symmetric trigonometric
polynomials dependent on the parameters
$q,t$. He came up with four main conjectures  formulated
for arbitrary root systems. A new approach to
the Macdonald theory was suggested in [C1] on the basis of
double affine Hecke algebras (new objects in mathematics).
In [C2] the norm conjecture
(including the famous constant term conjecture [M3])
and the conjecture
about the denominators of the coefficients of the Macdonald
polynomials were proved.
This paper contains the proof of
the remaining two (the duality and evaluation conjectures).

The evaluation conjecture (now a theorem) is in fact
a $q,t$-generalization
of the classic Weyl dimension formula. One can expect
interesting applications of this theorem since
the so-called $q$-dimensions are undoubtedly
important. It is  likely that we can incorporate
the Kac-Moody case as well. The necessary technique
was developed in [C4].

As to the  duality theorem (in its complete form), it states
that the generalized trigonometric-difference
zonal Fourier transform
is self-dual (at least formally).
We define this $q,t$-transform in terms of double affine
Hecke algebras. The most natural way to check the self-duality is to
use the connection of these algebras with the so-called
elliptic braid groups (the Fourier involution will turn into the
transposition of the periods of an elliptic curve).
\vskip 10pt

The classical trigonometric-differential Fourier transform
(corresponding to the limit $q=t^k$ as $t\to 1$ for certain special
$k$)
plays one of the main roles in the harmonic analysis
on symmetric spaces. It sends symmetric
trigonometric polynomials to the
corresponding radial  parts of
Laplace operators (Harish-Chandra,
Helgason) and is not self-dual. The calculation
of its inverse (the Plancherel theorem) is always
challenging and involving.

In the rational-differential setting, Charles Dunkl introduced
the generalized Hankel transform which appeared to be
 self-dual [D,J].
We demonstrate in this paper that one can save this very
important property if trigonometric polynomials
 come together with  difference operators.
At the moment, it is mostly an algebraic observation (the
difference-analitical aspects were not touched upon).

The root systems of type $A_n$ are always rather special.
First of all, we note the $q \leftrightarrow t$ symmetry and
very interesting positivity conjectures
(Macdonald [M1], Garsia, Haiman [GH]). Then the Macdonald
polynomials can be interpreted as generalized characters
(Etingof, Kirillov [EK1]). The difference Fourier transform
also has  particular features (we discuss this a little
at the end of the paper). By the way, due to Andrews
one can add $n$ new parameters
$q$
and still the
constant term conjecture (proved by Bressoud and
Zeilberger) will hold,  but there are no
related orthogonal polynomials.
In the differential setting, the
corresponding symmetric polynomials (Jack polynomials)
are  quite remarkable as well
(Hanlon, Stanley).

As to the differential theory, the Macdo\-nald- Mehta
conjectures
 were proved finally by Eric Opdam
[O1] (see also [O2]) excluding the duality con\-jec\-ture which
collapses
(the Fourier transform is not self-dual!). He used
the Heckman-Opdam operators
 (including the shift operator - see [01,He]).
We use their dif\-fe\-rence counterparts from [C1,C2]
defined by means of double
affine Hecke algebras. We mention that the latter algebras were
 not absolutely necessary in [C2] to prove the norm conjecture
 (the classic affine
Hecke algebras are enough). Only
in this paper the double Hecke algebras  work at their full potential
to ensure the duality.

\vskip 10pt
We note that this paper is a part of a new program
in the harmonic analysis of
symmetric spaces based  on certain
remarkable representations of  Hecke algebras
in terms of Dunkl and
Demazure operators instead of Lie groups and Lie algebras.
It gave already a parametric deformation of the
the classical  theory (see
[O1,He,C5]) directly connected with the so-called
quantun many-body problem
(Calogero,  Sutherland, Moser, Olshanetsky,
 Perelomov). Then it was extended
(in the algebraic context)
 to the difference, elliptic,
and finally to
the difference-elliptic case [C4]
presumably corresponding to the quantum Kac-Moody algebras.
Presumably because the harmonic analysis
for the latter algebras does not exist.

\vfil
{\it The duality-evaluation conjecture.}
Let $R=\{\al\}   \subset \R^n$ be a root system of type $A,B,...,F,G$
with respect to a euclidean form $(z,z')$ on $\R^n \ni z,z'$,
$W$ the Weyl group  generated by the the reflections $s_\al$.
We assume that $(\al,\al)=2$ for long $\al$.
Let us  fix the set $R_{+}$ of positive  roots ($R_-=-R_+$),
the corresponding simple
roots $\al_1,...,\al_n$, and  their dual counterparts
$a_1 ,..., a_n,  a_i =\al_i^\vee, \where \al^\vee =2\al/(\al,\al)$.
The dual fundamental weights
$b_1,...,b_n$  are determined from the relations  $ (b_i,\al_j)=
\de_i^j $ for the
Kronecker delta. We will also introduce the dual root system
$R^\vee =\{\al^\vee, \al\in R\}, R^\vee_+$, and the lattices
$$
\eqalignno{
& A=\oplus^n_{i=1}\Z a_i \subset B=\oplus^n_{i=1}\Z b_i,
}
$$
  $A_\pm, B_\pm$  for $\Z_{\pm}=\{m\in\Z, \pm m\ge 0\}$
instead of $\Z$. (In the standard notations, $A= Q^\vee,\
B = P^\vee $ - see [B].)  Later on,
$$
\eqalign{
&\nu_{\al}=\nu_{\al^\vee}\ =\ (\al,\al),\  \nu_i\ =\ \nu_{\al_i}, \
\nu_R\ = \{\nu_{\al}, \al\in R\}, \cr
&\rho_\nu\ =\ (1/2)\sum_{\nu_{\al}=\nu} \al \ =
\ (\nu/2)\sum_{\nu_i=\nu}  b_i, \for\al\in R_+,\cr
&r_\nu\ =\ \rho_\nu^\vee \ =\ (2/\nu)\rho_\nu\ =\
\sum_{\nu_i=\nu}  b_i,\quad 2/\nu=1,2,3.
}
\eqnu
$$

Let us put formally $x_i=exp({b_i}),\  x_b=exp(b)= \prod_{i=1}^n
x_i^{k_i} \for b=\sum_{i=1}^n k_i b_i$,  and introduce the algebra
$\C(\de,q)[x]$  of polynomials in terms of $x_i^{\pm 1}$ with the
coefficients belonging to the field $\C(\de,q)$ of rational functions
in terms  of indefinite complex parameters $\de, q_\nu,
\nu\in \nu_R$ (we will put $q_\al=q_{\nu_\al}=q_{\al^\vee}$).
The coefficient of $x^0=1$ ({\it the constant term})
will be denoted by $\langle \  \rangle$. The following product is a
Laurent series in $x$ with the coefficients in  $\C(\de,q)$:
$$
\eqalign{
&\mu\ =\ \prod_{a \in R_+^\vee}
\prod_{i=0}^\infty {(1-x_a\de_a^{i}) (1-x_a^{-1}\de_a^{i+1})
\over
(1-x_a q_a\de_a^{i}) (1-x_a^{-1}q_a^{-1}\de_a^{i+1})},
}
\eqnu
\label\mu\eqnum*
$$
where $\de_a=\de_{\nu}=\de^{2/\nu} \for \nu=\nu_a$.
We note that  $\mu\in
\C(\de,q)[x]$ if $q_\nu=\de_\nu^{k_\nu}$ for $k_\nu\in \Z_+$.

The
{\it monomial symmetric polynomials}
$m_{b}\ =\ \sum_{c\in W(b)}x_{c}$ for $b\in B_-$
form a base of the space
 $\C[x]^W$ of all $W$-invariant polynomials.
Setting  $\bar{x}_b\equal x_{-b}$,
$$
\eqalignno{
&\langle f,g\rangle\ =\langle \mu f\ \bar{g}\rangle \for
f,g \in \C(\de,q)[x]^W,
&\eqnu
}
$$
 we  introduce the {\it Macdonald
polynomials} $p_b(x),\   b \in B_-$, by means of
the conditions
$$
\eqalignno{
&p_b-m_b\ \in\ \oplus_c\C(\de,q)m_{c},\
\langle p_b, m_{c}\rangle = 0, &\eqnu
 \cr
&\where
 c\in B_-,\  c-b \in A_+, c\neq b.
\label\macd\eqnum*
}
$$
They can be determined by the Gram - Schmidt process
because the  pairing (see [M1,M2])  is non-degenerate
 and form a
basis in $\C(\de,q)[x]^W$.
Let $x_{i}(q^{-\rho}\de^{b})=
\de^{(b,b_i)}\prod_\nu q_\nu^{-(b_i,\rho_\nu)}$.

\proclaim{Main Theorem}
Given $b,c\in B_-$ and the corresponding Macdonald
polynomials $p_b, p_c$,
$$
\eqalignno{
&p_b(q^{-\rho} \de^{c})p_c(q^{-\rho})\ =\
p_c(q^{-\rho} \de^{b})p_b(q^{-\rho}),
&\eqnu\cr
\label\PP\eqnum*
&p_b(q^{-\rho})\ =\
\prod_\nu q_\nu^{(\rho_\nu,b)}
\prod_{a\in R_+^\vee, 0\le j< \infty} &\eqnu \cr
&\Bigl(
{
(1-\de_a^{j-(b,a^\vee)}\prod_\nu q_\nu^{(\rho_\nu,a)})
(1-q_a\de_a^{j}\prod_\nu q_\nu^{(\rho_\nu,a)})
 \over
(1-q_a\de_a^{j-(b,a^\vee)}\prod_\nu q_\nu^{(\rho_\nu,a)})
(1-\de_a^{j}\prod_\nu q_\nu^{(\rho_\nu,a)})
}
\Bigr).
\label\EV\eqnum*
}
$$
\label\MAIN\theoremnum*
\endproclaim
%\proofbox
The right hand side of (\ref\EV) is a rational function
in terms of $\de,q$ (we used $a^\vee=2a/(a,a)$ to make it more
transparent). We mention that there is a straightforward passage
 to the case where $\mu$ is introduced for
$\al\in R_+$ instead of $a\in R_+^\vee$ (see [C2]) and to
non-reduced root systems.

The second formula was conjectured by Macdonald (see (12.10),[M2]).
 He also formulated an equivalent version of (\ref\PP)
in one of his lectures (1991). Both statements seem to be
 esatblished in 1988 by Koornwinder for $A_n$ (his proof
was not published) and by Macdonald (to be published).
Recently the paper by Etingof and Kirillov [EK2]
appeared  were they use their interpretation of the Macdonald
polynomials to check the above theorem (and the norm conjecture)
in the case of $A_n$. As to other root systems,
it seems that almost
nothing was
known (excluding $BC_1$ and certain
 special values of the
parameters).

The author thanks G. Heckman
and A. Kirillov, Jr. for useful discussion.

%
%
%		Section 1
%
%
%\vskip 10pt
\section { Double affine Hecke algebras}
The vectors $\ \tal=[\al,k] \in
\R^n\times \R \subset \R^{n+1}$
for $\al \in R, k \in \Z $
form the {\it affine root system}
$R^a \supset R$ ( $z\in \R^n$ are identified with $ [z,0]$).
We add  $\al_0 \equal [-\th,1]$ to the  simple roots
for the {\it maximal root} $\th \in R$.
The corresponding set $R^a_+$ of positive roots coincides
with $R_+\cup \{[\al,k],\  \al\in R, \  k > 0\}$.

We denote the Dynkin diagram and its affine completion with
$\{\al_j,0 \le j \le n\}$ as the vertices by $\Ga$ and $\Ga^a$.
Let $m_{ij}=2,3,4,6$\  if $\al_i\and\al_j$ are joined by 0,1,2,3 laces
respectively.
The set of
the indices of the images of $\al_0$ by all
the automorphisms of $\Ga^a$ will be denoted by $O$ ($O=\{0\}
\for E_8,F_4,G_2$). Let $O^*={r\in O, r\neq 0}$.
The elements $b_r$ for $r\in O^*$ are the so-called minuscule
weights ($(b_r,\al)\le 1$ for
$\al \in R_+$).

Given $\tal=[\al,k]\in R^a,  \ b \in B$, let
$$
\eqalignno{
&s_{\tal}(\tz)\ =\  \tz-(z,\al^\vee)\tal,\
\ b'(\tz)\ =\ [z,\ze-(z,b)]
&\eqnu
%&(1.1)
}
$$
for $\tz=[z,\ze] \in \R^{n+1}$.

The {\it affine Weyl group} $W^a$ is generated by all $s_{\tal}$
(we write $W^a = <s_{\tal}, \tal\in R_+^a>)$. One can take
the simple reflections $s_j=s_{\al_j}, 0 \le j \le n,$ as its
generators and introduce the corresponding notion of the
length. This group is
the semi-direct product $W\lsmash A'$ of
its subgroups $W=<s_\al,
\al \in R_+>$ and $A'=\{a', a\in A\}$, where
$$
\eqalignno{
& a'=\ s_{\al}s_{[\al,1]}=\ s_{[-\al,1]}s_{\al}\for a=\al^{\vee},
\ \al\in R.
&\eqnu
%&(1.2)
}
$$

The {\it extended Weyl group} $ W^b$ generated by $W\and B'$
(instead of $A'$) is isomorphic to $W\lsmash B'$:
$$
\eqalignno{
&(wb')([z,\ze])\ =\ [w(z),\ze-(z,b)] \for w\in W, b\in B.
%&(1.3)
&\eqnu
}
$$

 Given $b_+\in B_+$, let
$$
\eqalignno{
&\om_{b_+} = w_0w^+_0  \in  W,\ \pi_{b_+} =
b'_+(\om_{b_+})^{-1}
\ \in \ W^b, \ \om_i=\om_{b_i},\pi_i=\pi_{b_i},
&\eqnu
\label\w0\eqnum*
}
$$
where $w_0$ (respectively, $w^+_0$) is the longest element in $W$
(respectively, in $ W_{b_+}$ generated by $s_i$ preserving $b_+$)
relative to the
set of generators $\{s_i\}$ for $i >0$.

We will  use here only the
elements $\pi_r=\pi_{b_r}, r \in O$. They leave $\Ga^a$ invariant
and form a group denoted by $\Pi$,
 which is isomorphic to $B/A$ by the natural
projection $\{b_r \to \pi_r\}$. As to $\{\om_r\}$,
they preserve the set $\{-\th,\al_i, i>0\}$.
The relations $\pi_r(\al_0)= \al_r= (\om_r)^{-1}(-\th)
$ distinguish the
indices $r \in O^*$. Moreover (see e.g. [C2]):
$$
\eqalignno{
& W^b  = \Pi \lsmash W^a, \where
  \pi_rs_i\pi_r^{-1}  =  s_j \if \pi_r(\al_i)=\al_j,\  0\le j\le n.
&\eqnu
%&(1.6)
}
$$

We extend the notion
of the length to $W^b$.
Given $\nu\in\nu_R,\  r\in O^*,\  \tw \in W^a$, and a reduced
decomposition $\tw\ =\ s_{j_l}...s_{j_2} s_{j_1} $ with respect to
$\{s_j, 0\le j\le n\}$, we call $l\ =\ l(\hw)$ the {\it length} of
$\hw = \pi_r\tw \in W^b$. Setting
$$
\eqalign{
\la(\hw) = &\{ \tal^1=\al_{j_1},\
\tal^2=s_{j_1}(\al_{j_2}),\
\tal^3=s_{j_1}s_{j_2}(\al_{j_3}),\ldots \cr
&\ldots,\tal^l=\tw^{-1}s_{j_l}(\al_{j_l}) \},
}
\eqnu
$$
\label\tal\eqnum*
one can represent
$$
\eqalign
{
&l=|\la(\hw)|=\sum_\nu l_\nu, \for l_\nu = l_\nu(\hw)=|\la_\nu(\hw)|,\cr
&\la_\nu(\hw) = \{\tal^{m},\ \nu(\tal^{m})= \nu(\tal_{j_m})= \nu\},
1\le m\le l,
}
\eqnu
\label\laset\eqnum*
$$
where $|\ |$  denotes the  number of elements,
 $\nu([\al,k]) \equal \nu_{\al}$.

For instance,
$$
\eqalign{
&l_\nu(b')\ =\ \sum_{\al} |(b,\al)|,\  \al\in R_+,
\nu_\al=\nu \in \nu_R, \cr
&l_\nu(b_+')\ =\ 2(b_+,\rho_\nu) \when b_+ \in B_+.
}
\eqnu
$$
\label\lb\eqnum*
Here $|\ | = $
absolute value. Later on  $b$ and $b'$ will not be distinguished.

We put
$m=2 \for D_{2k} \and C_{2k+1},\ m=1 \for C_{2k}, B_{k}$,
otherwise $m=|\Pi|$. The definition involves
the parameters $\de,\{ q_\nu , \nu \in \nu_R \}$ and independent
variables $X_1,\ldots,X_n$.
Let us set
$$
\eqalignno{
&   q_{\tal} = q_{\nu(\tal)},\ q_j = q_{\al_j},
\where \tal \in R^a, 0\le j\le n, \cr
& X_{\tb}\ =\ \prod_{i=1}^nX_i^{k_i}\de^{ k}
\if \tb=[b,k],
&\eqnu \cr
&\for b=\sum_{i=1}^nk_i b_i\in B,\ k \in {1\over m}\Z.
}
$$
\label\Xde\eqnum*

 Later on $ \C_{\de}$
 is the field of rational
functions in $\de^{1/m},$
$\C_\de[X] = \C_\de[X_b]$  means the algebra of
polynomials in terms of $X_i^{\pm 1}$
with the coefficients depending
on $\de^{1/m}$ rationally. We replace $ \C_{\de}$
by $ \C_{\de,q}$ if the functions (coefficients)
also depend rationally
on $\{q_\nu^{1/2} \}$.

Let $([a,k],[b,l])=(a,b)$ for $a,b\in B,\
[\al,k]^\vee= [\al^\vee,k],\ a_0=\al_0,\ \nu_{\al^\vee}=\nu_\al, $
and  $ \al_{r^*} \equal \pi_r^{-1}(\al_0)$ for $  r\in O^*$.

\proclaim{Definition }
 The  double  affine Hecke algebra $\HH\ $
(see [C1,C2])
is generated over the field $ \C_{\de,q}$ by
the elements $\{ T_j,\ 0\le j\le n\}$,
pairwise commutative $\{X_b, \ b\in B\}$ satisfying (\ref\Xde),
 and the group $\Pi$ where the following relations are imposed:

(o)\ \  $ (T_j-q_j^{1/2})(T_j+q_j^{-1/2})\ =\ 0,\ 0\ \le\ j\ \le\ n$;

(i)\ \ \ $ T_iT_jT_i...\ =\ T_jT_iT_j...,\ m_{ij}$ factors on each side;

(ii)\ \   $ \pi_rT_i\pi_r^{-1}\ =\ T_j \if \pi_r(\al_i)=\al_j$;

(iii)\  $T_iX_b T_i\ =\ X_b X_{a_i}^{-1} \if (b,\al_i)=1,\
1 \le i\le  n$;

(iv)\  $T_0X_b T_0\ =\ X_{s_0(b)}\ =\ X_b X_{\th}\de^{-1}
\if (b,\th)=-1$;

(v)\ \ $T_iX_b\ =\ X_b T_i$ if $(b,\al_i)=0 \for 0 \le i\le  n$;

(vi)\ $\pi_rX_b \pi_r^{-1}\ =\ X_{\pi_r(b)}\ =\ X_{\om^{-1}_r(b)}
\de^{(b_{r^*},b)},\  r\in O^*$.
\label\double\theoremnum*
\endproclaim
\proofbox

Given $\tw \in W^a, r\in O,\ $ the product
$$
\eqalignno{
&T_{\pi_r\tw}\equal \pi_r\prod_{k=1}^l T_{i_k},\where
\tw=\prod_{k=1}^l s_{i_k},
l=l(\tw),
&\eqnu
\label\Tw\eqnum*
}
$$
does not depend on the choice of the reduced decomposition
(because $\{T\}$ satisfy the same ``braid'' relations as $\{s\}$ do).
Moreover,
$$
\eqalignno{
&T_{\hv}T_{\hw}\ =\ T_{\hv\hw}\  \hbox{ whenever}\
 l(\hv\hw)=l(\hv)+l(\hw) \for
\hv,\hw \in W^b.
 %&(2.7)}
&\eqnu}
$$
\label\TT\eqnum*
  In particular, we arrive at the pairwise
commutative elements
$$
\eqalignno{
& Y_{b}\ =\  \prod_{i=1}^nY_i^{k_i} \if
b=\sum_{i=1}^nk_ib_i\in B,\where
 Y_i\equal T_{b_i},
&\eqnu
\label\Yb\eqnum*
}
$$
satisfying the relations
$$
\eqalign{
&T^{-1}_iY_b T^{-1}_i\ =\ Y_b Y_{a_i}^{-1} \if (b,\al_i)=1,
\cr
& T_iY_b\ =\ Y_b T_i \if (b,\al_i)=0, \ 1 \le i\le  n.}
%\eqno(2.9)
\eqnu
$$
Let us introduce the following elements from
$\C_q^n$:
$$
\eqalign{
&q^{\pm\rho}\equal (l_q(b_1)^{\pm 1},\ldots,l_q(b_n)^{\pm 1}),\where\cr
&l_q(\hw)\equal \ \prod_{\nu\in\nu_R} q_\nu^{l_\nu(\hw)/2},\
\hw\in W^b,
}
\eqnu
\label\qlen\eqnum*
$$
and the corresponding {\it evaluation maps}:
$$
\eqalign{
&X_i(q^{\pm\rho})= l_q(b_i)^{\pm 1} = Y_i(q^{\pm\rho}),\ 1\le i\le n.
}
\eqnu
\label\eval\eqnum*
$$
For instance, $X_{a_i}(q^{\rho})\ =\ l_q(a_i)= q_i$ (see (\ref\lb)).

\proclaim {Theorem}
i) The elements $H \in \HH\ $  have
the unique decompositions
$$
\eqalignno{
&H =\sum_{w\in W }  g_{w}  T_{w} f_w,\
g_{w} \in \C_{\de,q}[X],\ f_{w} \in \C_{\de,q}[Y].
&\eqnu
}
$$

ii) The   map
$$
\eqalign{
 \vph: &X_i \to Y_i^{-1},\ \  Y_i \to X_i^{-1},\  \ T_i \to T_i, \cr
&q_\nu \to q_\nu,\
\de\to \de,\ \nu\in \nu_R,\ 1\le i\le n.
}
\eqnu
\label\vph\eqnum*
$$
can be extended to an anti-involution
($\vph(AB)=\vph(B)\vph(A)$)
 of \HH\ .

iii) The  linear functional on \HH\
$$
\eqalignno{
&[\![ \sum_{w\in W }  g_{w}  T_{w} f_{w}]\!]\ =\
\sum_{w\in W} g_{w}(q^{-\rho}) l_q(w) f_{w}(q^{\rho})
&\eqnu
\label\brack\eqnum*
}
$$
is invariant with respect to $\vph$. The bilinear form
$$
\eqalignno{
&[\![ F,G]\!]\equal [\![ \vph(G)H]\!],\
G,H\in \HH\ ,
&\eqnu
\label\form\eqnum*
}
$$
is symmetric ($[\![ G,H]\!]= [\![ H,G]\!]$)
and non-degenerate.
\endproclaim
\label\dual\theoremnum*

{\it Proof.} The first statement is  from Theorem 2.3 [C2].
The map $\vph$ is the composition of the
involution (see [C1])
$$
\eqalign{
  &X_i \to Y_i,\ \  Y_i \to X_i,\  \ T_i \to T_i^{-1}, \cr
&q_\nu \to q_\nu^{-1},\
\de\to \de^{-1},\ 1\le i\le n,
}
\eqnu
$$
and the main anti-involution "*" from [C2], sending
$$
\eqalign{
  & X_i \to X_i^{-1},\ \  Y_i \to Y_i^{-1},\  \ T_i \to T_i^{-1}, \cr
&q_\nu \to q_\nu^{-1},\
\de\to \de^{-1},\ 0\le i\le n.
}
\eqnu
$$
The other claims follow directly from the definition of
$[\![ \ ]\!]$.
\proofbox

One can extend $[\![\ ]\!]$ to the localization of \HH\
with respect to all polynomials in $X$ (or in $Y$). The
algebra becomes the semi-direct product of $\C[W^b]$ and
$\C(X)$ after this (see [C3]). Sometimes it is also convenient
to involve proper completions of $\C(X)$ (see the end of the
paper).

%
%
%		Section 2
%
%
%\vskip 10pt
\section { Difference operators}
Setting (see the Introduction)
$$
\eqalignno{
& x_{\tb}=  \prod_{i=1}^nx_i^{k_i}\de^{ k} \if
\tb=[b,k],
b=\sum_{i=1}^nk_i b_i\in B,\ k \in {1\over m}\Z,
&\eqnu
\label\xde\eqnum*}
$$
for independent $x_1,\ldots,x_n$, we will
 consider $\{X\}$ as  operators acting in $\C_\de[x]=
\C_\de[x_1^{\pm 1},\ldots,x_n^{\pm 1}]$:
$$
\eqalignno{
& X_{\tb}(p(x))\ =\ x_{\tb} p(x),\    p(x) \in
\C_\de [x].
&\eqnu}
$$
\label\X\eqnum*
The elements $\hw \in W^b$ act in $\C_{\de}[x]$
 by the
formulas:
$$
\eqalignno{
&\hw(x_{\tb})\ =\ x_{\hw(\tb)}.
&\eqnu}
$$
 In particular:
$$
\eqalignno{
&\pi_r(x_{b})\ =\ x_{\om^{-1}_r(b)}\de^{(b_{r^*},b)}
\for \al_{r^*}\ =\ \pi_r^{-1}(\al_0), \ r\in O^*.
&\eqnu}
$$
\label\pi\eqnum*

The {\it Demazure-Lusztig operators} (see
[KL, KK, C1], and [C2] for more detail )
$$
\eqalignno{
&\hT_j\  = \  q_j ^{1/2} s_j\ +\
(q_j^{1/2}-q_j^{-1/2})(X_{a_j}-1)^{-1}(s_j-1),
\ 0\le j\le n.
&\eqnu
\label\Demaz\eqnum*
}
$$
act   in $\C_{\de,q}[x]$ naturally.
We note that only $\hT_0$ depends on $\de$:
$$
\eqalign{
&\hT_0\  =  q_0^{1/2}s_0\ +\ (q_0^{1/2}-q_0^{-1/2})
(\de X_{\th}^{-1} -1)^{-1}(s_0-1),\cr
&\where
s_0(X_i)\ =\ X_iX_{\th}^{-(b_i,\th)}\de^{(b_i,\th)}.
}
%\eqno(2.12)
\eqnu
$$

\proclaim{Theorem }
 The map $ T_j\to \hT_j,\ X_b \to X_b$ (see (\ref\Xde,\ref\X)),
$\pi_r\to \pi_r$  (see (\ref\pi)) induces a $ \C_{\de,q}$-linear
homomorphism from \HH\ to the algebra of linear endomorphisms
of $\C_{\de,q}[x]$.
 This representation is faithful and
remains faithful when   $ \de,q$ take  any non-zero
values assuming that
 $\de$ is not a root of unity (see [C2]). The image $\hat{H}$
is uniquely determined from the following condition:
$$
\eqalign{
&\hat{H}(f(x))\ =\ g(x)\for H\in \HH\ ,\if Hf(X)\ =\cr
 &g(X)+ \sum_{i=0}^n H_i(T_i-q_i)+
\sum_{r\in O^*} H_r(\pi_r-1), \where H_i,H_r\in \HH\ .
}
\eqnu
\label\hat\eqnum*
$$
\endproclaim
\proofbox
\label\faith\theoremnum*

Due to  Theorem \ref\dual, an arbitrary  $H\in \HH\ $ can be
uniquely represented in the form
$$
\eqalign{
H =&\sum_{b\in B, w\in W }  g_{b,w}  Y_b T_{w},\
g_{b,w} \in \C_{\de,q}[X],\cr
=&\sum_{b\in B, w\in W }   T_{w} X_b  g'_{b,w},\
g'_{b,w} \in \C_{\de,q}[Y].
}
\eqnu
$$
We set:
$$
\eqalign{
&[H]_{\dagger}\  =\ \sum_{b\in B, w\in W }  g_{b,w}  Y_b l_q(w),\
{}_{\dagger}[H]\ =\ \sum_{b\in B, w\in W }   l_q(w) X_b  g_{b,w},\cr
&[H]_\ddagger\  =\ \sum_{b\in B, w\in W }  g_{b,w}  [\![Y_bT_w]\!],\
{}_\ddagger[H]\ =\ \sum_{b\in B, w\in W }    [\![T_wX_b]\!]  g'_{b,w}.
}
\eqnu
\label\redH\eqnum*
$$
One easily checks that
$$
\eqalign{
[\![H_1 H_2]\!]\ =\ &[\![H_1 [H_2]_{\dagger}]\!]\ =\
[\![ {}_{\dagger}[H_1]H_2]\!]\ =\cr
&[\![H_1 [H_2]_\ddagger]\!]\ =\ [\![ {}_\ddagger[H_1]H_2]\!]
\for H_1,H_2\in \HH\ .
}
\eqnu
\label\dag\eqnum*
$$

Let us represent the image $\hH$ of $H$ as follows:
$$
\eqalignno{
&\hH\ = \sum_{b\in B, w\in W} h_{b,w} b  w,
\ =\ \sum_{b\in B, w\in W} w b h'_{b,w}.
&\eqnu
\label\hatH\eqnum*
}
$$
where  $h_{b,w}, h'_{b,w}$ belong to  the field $\C_{\de,q}(X)$
of rational  functions in
$X_1,...,X_n$. We extend the above operations to arbitrary operators
in the form (\ref\hatH):
$$
\eqalign{
&[\hH]_{\dagger}= \sum h_{b,w}  b ,\
{}_{\dagger}[\hH] = \sum b h'_{b,w}\  ,\
[\![ \hH ]\!]\
 = \sum h_{b,w}
(q^{-\rho}).
}
\eqnu
\label\Brack\eqnum*
$$
These operations commute with the homomorphism $H\to\hH$.

Let us define the
{\it difference Harish-Chandra map} (see
[C2], Proposition 3.1):
$$
\eqalignno{
&\chi( \sum_{w\in W,b\in B}  h_{b,w} b w)\ =\
\sum_{b\in B,w\in B} h_{b,w}(\diamondsuit)
y_b \ \in \C_{\de,q}[y],
&\eqnu
}
$$
where
$\diamondsuit \equal ( X_1=...=X_n=0 ),\ \{y_b\}$ is one more  set of variables
 introduced  for
independent $y_1,...,y_n$ in the same way as
$\{x_b\}$ were.

\proclaim{Proposition}
 Setting
$$
\eqalign{
&\l_f \ =\ f(Y), \ \hat{\l}_f=\ f(\hY),\
L_f = L_f^{\de,q}\ \equal\ [(\hat{\l}_f)]_{\dagger}
}
\eqnu
\label\Ll\eqnum*
$$
for $f=\sum_b g_b y_b \in \C_{\de,q}[y]$, one has:
$$
\eqalign{
\chi(\hat{\l}_f)\ =\  \chi(L_f)\ =\
[\![ f(Y)]\!]\ =\
\sum_{b\in B}  g_{b}
\prod_\nu q_\nu^{(b,\rho_\nu)} y_b.
}
\eqnu
$$
\label\chi\eqnum*
\endproclaim
\proofbox

The proof
of the following theorem repeats the proof of Theorem 4.5,[C2]
(where the relations $q_\nu=\de_\nu^{k_\nu}$ for $k_\nu\in \Z_+$
were imposed). We note that once  (\ref\Lf) is known for these
special $q$ it holds true for all $\de,q$ since all the
coefficients of difference operators and polynomials
are rational in $\de,q$.

\proclaim{Theorem}
The difference operators $\{ L_f, \
f(y_1,\ldots,y_n)\in \C_{\de,q}[y]^W\}$
are pairwise commutative,  $W$-invariant (i.e $w L_f w^{-1}=$
$L_f$ for all $w\in W$) and preserve $\C_{\de,q}[x]^W$. The
Macdonald polynomials $p_b=p_b^{\de,q} (b\in B_-)$
from  (\ref\macd) are their eigenvectors:
$$
\eqalignno{
&L_f(p_b^{\de,q})=f(q^\rho\de^{-b}) p_b^{\de,q},\
y_i(q^\rho\de^{-b})\equal
\de^{-(b_i,b)}\prod_\nu q_\nu^{(b_i,\rho_\nu)}.
&\eqnu
\label\Lf\eqnum*}
$$
\label\LF\theoremnum*
\endproclaim
\proofbox

We fix a subset $v\in \nu_R$ and introduce the
{\it shift operator} by the formula
$$
\eqalignno{
&\g_v \ =\
(\x_v) ^{-1}\y_v,\ G_v^{\de,q}\ =\ [\hat{\g}_v]_{\dagger}  \ =
\ (\x_v )^{-1}[\hat{\y}_v]_{\dagger},
 &\eqnu
\label\shift\eqnum*
}
$$
$$
 \x_v  = \prod_{\nu_a\in v}((q_a X_{a})^{1/2}-
(q_a X_{a})^{-1/2}),\  \y_v  = \prod_{\nu_a\in v}
(q_a Y_{a}^{-1})^{1/2}-
(q_a Y_{a}^{-1})^{-1/2}).
$$
Here $a=\al^\vee\in R_+^\vee, \nu_{a}=\nu_\al, q_a=q_\al$, the
elements $\x_{v}= \x_{v}^q, \y_{v}=\y_{v}^q $
belong to $\C_q [X],
\C_q [Y]$ respectively.

\proclaim{Theorem}
The operators $\hat{\g}_{v}$
and $G_v^{\de,q}$  are $W$-inva\-riant and preserve $ \C_{\de,q}[x]^W$
(their restrictions to the latter space coincide). Moreover,
if $q_{\nu}=1$ when
$\nu\not\in v$ then
$$
\eqalign{
& G_v^{\de,q}L_f^{\de,q}\ =\ L_f^{\de, q\de_v} G_v^{\de,q}
\for f\in \C_{\de,q}[y]^W,
\cr
&G_v^{\de,q} (p_{b}^{\de,q})= g_v^{\de,q}(b)
p_{b+r_v}^{\de, q\de_v}, \for\cr
&g_v^{\de,q}(b)\ =\
\prod_{a\in R_+^\vee,\nu_a\in v} (y_a(q^{\rho/2}\de^{-b/2}) -
q_a y_a(q^{-\rho/2}\de^{+b/2})),
}
\eqnu
\label\Gnu\eqnum*
$$
where $r_v=\sum_{\nu\in v}r_\nu,\ q\de_v=\{q_\nu\de^{2/\nu} , q_{\nu'}\}$
for $\nu\in v\not\ni \nu'\ $, $p_c=0 \for c\not\in B_-$.
\label\Gp\theoremnum*
\endproclaim

{\it Proof.}  When $q_\nu=\de^{2k_\nu/\nu}$ for $k_\nu\in \Z_+$
these statements  are in fact from [C2].
They give (\ref\Gnu) for all $\de,q$. Indeed,
 it can be rewritten as follows:
$$
\eqalign{
& [\hat{\l}_f^{\de,q}\x_v^{q}]_\dagger\  =\
\x_v^{q} L_f^{\de,q\de_v},
}
\eqnu
\label\xGx\eqnum*
$$
where the coefficients of the difference operators on both
sides are from
$\C_{\de,q}[X]$.
Here we used that $[\l\m]_{\dagger}=[\l]_{\dagger}[\m]_{\dagger} $ for
arbitrary operators  $\l,\m$ in the form (\ref\hatH) if the
second is $W$-invariant. The remaining formulas can be
deduced from [C2] in the same way (they mean certain
identities in $\C_{\de,q}$ which are enough to check for
$q_\nu=\de^{2k_\nu/\nu}$). One can use (\ref\chi) as well.
\proofbox

%
%
%		Section 3
%
%
%\vskip 10pt
\section { Duality and evaluation conjectures}
First of all we will use Theorem \ref\dual to define the
{\it zonal Fourier transform}. We will sometimes
identify the elements
$H\in \HH\ $ with their images $\hH$. The following pairing
on $f,g\in \C_{\de,q}[x]$ is symmetric and non-degenerate:
$$
\eqalign{
 &[\![f,g]\!]= [\![f(X),g(X)]\!] = [\![\vph(f(X))g(X)]\!]  =\cr
&[\![\bar{f}(Y)g(X)]\!] =  \{\l_{\bar{f}}(g(x))\}(q^{-\rho}).
}
\eqnu
\label\Fourier\eqnum*
$$
Here $\bar{x}_b=x_{-b}=x_b^{-1}, \ {\l}$ is from (\ref\Ll),
and we used the main defining property (\ref\hat)
of the representation
from  Theorem \ref\faith. The pairing remains non-degenerate
when restricted to $W$-invariant polynomials.

\proclaim{ Definition}
The Fourier transforms $\vph(\l),\vph(L)$ of $\C_{\de,q}$-linear
operators acting respectively either in $\C_{\de,q}[x]$ or
in $\C_{\de,q}[x]^W$ are defined from the relations:
$$
\eqalign{
&[\![\l(f),g]\!]\ =\  [\![f,\vph(\l)(g)]\!],\ f,g\in C_{\de,q}[x],\cr
&[\![L(f),g]\!]\  =\   [\![f,\vph(L)(g)]\!],\ f,g\in \C_{\de,q}[x]^W.
}
\eqnu
$$
If $\l$ preserves $\C_{\de,q}[x]^W$ then so does $\vph(\l)$
and $\vph(L)=
[\vph(\l)]_{\dagger}, \where L=[\l]_{\dagger}$ is the restriction of
$\l$ to the invariant polynomials.
\endproclaim
\label\FT\theoremnum*
\proofbox

This involution  ($\vph^2=\hbox{id}$) extends  $\vph$ from
(\ref\vph) by construction. If $f\in C_{\de,q}[x]^W$, then
$\vph(L_f)= [\bar{f}(X)]_{\dagger}$. We arrive at the following
theorem:

\proclaim {Duality Theorem}
Given $b,c\in B_-$ and the corresponding Macdonald's
polynomials $p_b, p_c$,
$$
p_b(q^{-\rho} \de^{c})p_c(q^{-\rho})\ =
[\![p_b,p_c]\!]\ =\ [\![p_c,p_b]\!]\ = \
p_c(q^{-\rho} \de^{b})p_b(q^{-\rho}).
\eqnu
\label\pp\eqnum*
$$
\endproclaim
\proofbox

To complete this theorem we need to calculate $p_b(q^{-\rho})$.
The main step is the formula for
$p'((q\de_v)^{-\rho})$ in terms of $p(q^{-\rho})$,
where  (see (\ref\Gnu))
$$p=p_b,\ p'\ =\ p_{b+r_v}^{q\de_v},\
p'=(g_v^{q}(b))^{-1}
G_v^{q} (p).$$
Here and in similar formulas we show the dependence on $q$
omitting $\de$
since the latter  will be the same for all polynomials and
operators. Let
$$ \bar{\y}_v^q  = \prod_{a\in R_+^\vee,\nu_a\in v}
((q_a Y_{a})^{1/2}-
(q_a Y_{a})^{-1/2}).$$

\proclaim { Key Lemma}
$$
\eqalign{
&d_v^q p'((q\de_v)^{-\rho})\ = \cr
&\prod_{a\in R_+^\vee,\nu_a\in v} \Bigl(q_a^{-1} y_a(q^{-\rho/2}\de^{+b/2})
-y_a(q^{+\rho/2}\de^{-b/2})\Bigr)p(q^{-\rho}), \cr
&d_v^q\ =\
\prod_{a\in R_+^\vee,\nu_a\in v} \Bigl(q_a^{-1} y_a((q\de_v)^
{-\rho/2})
-y_a((q\de_v)^{+\rho/2})\Bigr)m_{-r_v}(q^{-\rho}).
}
\eqnu
\label\key\eqnum*
$$
\label\KEY\eqnum*
\endproclaim
{\it Proof.}
Let us use formula (\ref\xGx):
$$
\eqalign{
&[\![({\y}_v^q\x_v^q)(\x_v^q)^{-1}\l^q_{\bar{p}'}\x_v^q]\!]\ =
\ [\![({\y}_v^q\x_v^q)[(\x_v^q)^{-1}\l^q_{\bar{p}'}\x_v^q]_{\dagger}
]\!]\ =\cr
&[\![(\hat{\y}_v^q\hat{\x}_v^q)   \bar{p}'(\hat{Y}^{q\de_v})]\!]\ =
\ [\![({\y}_v^q\x_v^q) ]\!] p'((q\de_v)^{-\rho}).
}
\eqnu
\label\one\eqnum*
$$
On the other hand, it equals:
$$
\eqalign{
&[\![{\y}_v^q p'(Y) \x_v^q ]\!]\ =
[\![\y_v^q {p}'({X}^{q}) {\x}_v^q ]\!]\ =\cr
&[\![{\y}_v^q (\x_v^q p'(x)) ]\!]\ =
\pm[\![\bar{\y}_v^q (\x_v^q p'(x)) ]\!].
}
\eqnu
$$
Here we applied the anti-involution $\vph\  (\vph(\x)=\y,\
\vph(\y)=\x)$, then went
 from the abstract
$[\![\ ]\!]$ to that from (\ref\Brack), and
used  Theorem \ref\faith. The last transformation requires
 special comment. We will justify it in a moment.

After this, one can use  (\ref\Lf):
$$
\eqalign{
&[\![\bar{\y}_v^q (\x_v^q p'(x)) ]\!]\ =
[\![(\bar{\y}_v^q {\y}_v^q)   g_v^q(b)^{-1}p(x)) ]\!]\ =\cr
& g_v^q(b)^{-1}(\bar{\y}\y)(q^\rho\de^{-b})
[\![ p(x)]\!]\ =\cr
&\prod_{a\in R_+^\vee,\nu_a\in v} (q_a^{-1} y_a(q^{-\rho/2}\de^{+b/2}) -
y_a(q^{+\rho/2}\de^{-b/2}))[\![ p(x) ]\!].
%Y change by Y^{-1}, then q_a move with change by q_a^{-1} because
%the term from \x will be ()^{-1}
}
\eqnu
\label\evalp\eqnum*
$$
Finally, ${d}_v^q\equal \pm[\![{\y}_v^q\x_v^q ]\!]$ can be
determined from (\ref\evalp) and the relation
$1=p'\ =\ g_v^{q}(b)^{-1}
G_v^{q} (p_{b}^{q})\for b=-r_v,$ where $p_{-r_v}$ coincides with
the monomial function $m_{-r_v}$ (it follows directly from
the definition):
$$
\eqalign{
&{d}_v^q\ =\
\prod_{a\in R_+^\vee,\nu_a\in v} (q_a^{-1} y_a((q\de_v)^
{-\rho/2})
-y_a((q\de_v)^{+\rho/2})) m_{-r_v}(q^{-\rho}).
}
\eqnu
\label\dterm\eqnum*
$$

Let us check that
$$[\![(\bar{\y}_v^q-l_\ep(w_0){\y}_v^q)
(\x_v^q p'(x)) ]\!]=0 \for \hbox {\ any\ }
p'\in \C[x],
$$
where $l_\ep(w_0)= \prod_\nu \ep_\nu^{l_\nu(w_0)}$,
$$
\ep =   \{\ep_\nu= -1 \if \nu\in v,\hbox{\ otherwise\ } \ep_\nu=1\},\
 \nu\in \nu_R.
$$
Following formula (4.18),[C2] we introduce
the
{\it $q$-symmetrizers}, setting
$$
\eqalign{
&\p_v^q\ =\ (\pi_v^q)^{-1}\sum_{w\in W}
\prod_\nu(\ep_\nu q_\nu^{1/2})^{\ep_\nu(l_\nu(w)-l_\nu(w_0))} T_w,
\cr
&\pi_v^q\ =\ \sum_{w\in W}
\prod_\nu (\ep_\nu q_\nu^{1/2})^{\ep_\nu(2l_\nu(w)-l_\nu(w_0))}.
}
\eqnu
$$
It results from Proposition 3.5 and  Corollary 4.7(ibidem) that
$$\p_v^q (\x_v^q p')=\x_v^q p',\  \hat{\p}_v^q\p_v^{q=0}=\hat{\p}_v^q, \if
 q_\nu=1 \for
\nu\not\in \nu_R. $$
Hence
$$
\eqalign{
& [\![(\bar{\y}_v^q -l_\ep(w_0){\y}_v^q) (\x_v^q p'(x)) ]\!] =
[\![(\bar{\y}_v^q -l_\ep(w_0){\y}_v^q)\p_v^q  (\x_v^q p'(x)) ]\!] =\cr
&[\![ (\y_v^q \bar{p}'(Y))\p_v^q (\bar{\x}_v^q-l_\ep(w_0){\x}_v^q) ]\!] =
[\![ (\y_v^q \bar{p}'(Y))
\{\hat{\p}_v^q\p_v^{q=0} (\bar{\x}_v^q-l_\ep(w_0){\x}_v^q)\}
 ]\!].
}
$$
The latter equals zero.
\proofbox

Let us take any set $k=\{k_{\nu_1}\ge k_{\nu_2}\}\in \Z_+$ and put
$$
\eqalign{
q(k)=\{\de^{2k_\nu/\nu}\},\   k\cdot r=\sum_\nu k_\nu r_\nu,\
p_b^{(k)}= p_b^{q(k)}.
}
\eqnu
$$
The remaining  part of the calculation is based on
 the following chain of the shift operators that will be applied
to $p^{(0)}_{b-k\cdot r}=m_{b-k\cdot r}$
one after another:
$$
\eqalignno{
&G_{\nu_R}^{(k-1)}G_{\nu_R}^{(k-2)}\cdots
G_{\nu_R}^{(k-s)}G_{\nu_1}^{(k-s-e)}
\cdots G_{\nu_1}^{(0)},
&\eqnu
}
$$
where
$ k_{\nu_1}=s+t,\ k_{\nu_2}=s,\
e=\{e_\nu\},\  e_{\nu_1}=1,\ e_{\nu_2}=0,\ k-s= te$,
the set $\{1,1\}$ is denoted by $1$.

Lemma \ref\KEY gives that for a certain $D^{(k)}$ (which does
not depend on $b$):
$$
\eqalignno{
D^{(k)}p_b^{(k)}(q(k)^{-\rho})&\ =
m_{b-k\cdot r}(1)\prod_{a\in R_+^\vee,\nu_a\in v(i)}^{0\le i<s+t}
 &\eqnu \cr
&\Bigl(q(i)_a\, y_a(q(i)^{\rho/2}\de^{-b(i)/2})
-y_a(q(i)^{-\rho/2}\de^{b(i)/2})\Bigr),
\label\prodb\eqnum*
}
$$
where $q(i)=q(k(i)),\ b(i)=b -(k-k(i))\cdot r$,
$$k(i)=ie,\ v(i)=\nu_1\if i < t, \
 k(i)= i-t + te, \
v(i)=\nu_R \if i\ge t.$$

As to $D^{(k)}$, it equals the right hand side of (\ref\prodb)
when $b=0$. We note that
$$q(i)_a= \de^{2j/\nu_a} \for j=k_a+i-s-t,\ k_a=k_{\nu_a},$$
because  $i\ge t\if \nu_a\neq \nu_1$ (and $0\le j< k_a$).
The relation  $(2/\nu)\rho_\nu=r_\nu$ leads to the formulas:
$$
\eqalign{
&q(i)^{\rho/2}\de^{-b(i)/2}\ =\
\de^{(k(i)\cdot r -b+(k-k(i))\cdot r)/2}\ =\  \de^{k\cdot r-b/2},\cr
&y_a(q(i)^{\rho/2}\de^{-b(i)/2})\ =\
\de^{ (k\cdot r-b,a)/2 }.
}
\eqnu
$$
 Finally, we arrive at the following theorem:

\proclaim{ Evaluation Theorem}
$$
\eqalignno{
&p_b^{(k)}(q(k)^{-\rho})=
{m_{b-k\cdot r}(1)\over m_{-k\cdot r}(1)}
\prod_{\al\in R_+, 0\le j< k_\al} &\eqnu \cr
&\Bigl(
{
\de^{ \{(k\cdot r-b,\al)+j\}/\nu_\al }
-\de^{- \{(k\cdot r-b,\al)+j\}/\nu_\al}
 \over
\de^{ \{(k\cdot r,\al)+j\}/\nu_\al }
-\de^{- \{(k\cdot r,\al)+j\}/\nu_\al}
}
\Bigr).
\label\evde\eqnum*
}
$$
\endproclaim
\proofbox

We note that
$m_{b-k\cdot r}(1)/ m_{-k\cdot r}(1)\ =\
|W(b-k\cdot r)|/|W(k\cdot r)|\in \Z_+$. It  equals $1$
for all $b\in B_-$ when
$\prod_\nu k_\nu\neq 0$. Assuming this
we have:
$$
\eqalignno{
&p_b^{q(k)}(q(k)^{-\rho})\ =\
\de^{(k\cdot r,b)}
\prod_{\al\in R_+, 0\le j< \infty} &\eqnu \cr
&\Bigl(
{
(1-\de^{ 2\{(k\cdot r-b,\al)+j\}/\nu_\al })
(1-q_\al(k)\de^{ 2\{(k\cdot r,\al)+j\}/\nu_\al })
 \over
(1-q_\al(k)\de^{ 2\{(k\cdot r-b,\al) +j\}/\nu_\al })
(1-\de^{ 2\{(k\cdot r,\al)+j\}/\nu_\al })
}
\Bigr).
\label\evdeq\eqnum*
}
$$
The limit of (\ref\evdeq) as one of the $k_\nu$ approaches zero
exists and coincides with (\ref\evde). Since both sides of
this formula are rational functions in $q(k)\and \de$ we get
 (\ref\EV)  (cf.\ Theorem \ref\Gp).

We note that actually this paper does not depend very much
on the definition of the Macdonald polynomials from the
Introduction. We can eliminate $\mu$ introducing these
polynomials as the eigenfunctions of the $L$-operators
(formula (\ref\Lf)). Therefore it is likely that
paper [C4] can be extended to give  a "difference-elliptic"
Weyl dimension formula.

\vskip 10pt
 {\it Schwartz functions.}
 In conclusion we will use Macdonald's polynomials
to construct  pairwise orthogonal functions
 with respect to the pairing
$[\![\ ,\ ]\!]$ in the case of $A_n$.
At the moment, the extension of this construction
to other root systems is not known. We will start with
the following observation:

\proclaim{ Proposition}
Adding  proper roots of $\de$, the following maps
are automorphisms of \HH\ of type $A_n$:
$$
\eqalign{
& \tau: \ X_i \to X_i,\ \ Y_i \to X_iY_i\de^{-c_i},\
\ T_i \to T_i, \cr
& \om: \ X_i \to Y_i,\ \ Y_i \to Y_i^{-1}X_i^{-1}Y_i\de^{2c_i},\
\ T_i \to T_i, \cr
&q_\nu \to q_\nu,\
\de\to \de,\ c_i= (b_i,b_i)= i(n-i+1)/(2(n+1)).
}
\eqnu
\label\tauom\eqnum*
$$
\endproclaim
{\it Proof.} Here $b_i=\om_i$ in the  notations from [B].
%c_i =(b_i,\rho_2)/(1+(\th,\rho_2)).
 The proof can be deduced
from the topological interpretation of \HH\ from [C1]
in terms of the elliptic braid groups. In the case of $A_n$
the latter group  (due to
 Birman and  Scott)
is especially simple. Actually these automorphisms
are related to the standard generators of
$SL_2(\Z)$. Let us give another description of $\tau$
(as for $\om$, it can be expressed in terms of $\tau$ and
$\vph$).

Setting
$x_b=\de^{z_b}, \ z_{a+b}=z_a+z_b, \ z_i=z_{b_i},\
a(z_b)= z_b-(a,b), \ a,b\in \R^n,$
we introduce the {\it Gaussian function}
$\ga\ =\ \de^{\Sigma_{i=1}^n z_i z_{\al_i}/2}$,
which is considered as a  formal series in $x, \log\de$
and satisfies the following difference relations:
$$
\eqalign{
&b_j(\ga)\ =\ \de^{(1/2)\Sigma_{i=1}^n (z_i-(b_j,b_i))
(z_{\al_i}-\de_i^j)}\ =\cr
& \ga \de^{-z_j+ (b_j,b_j)/2 }\ =\  x_j^{-1}\ga \de^{(b_j,b_j)/2}
\for 1\le j\le n.
}
\eqnu
\label\gauss\eqnum*
$$
The Gaussian function commutes with $T_j,\for 1\le j\le n$
because it is $W$-invariant.
Since all $b_j$ are minuscule, we  use directly  formulas
(\ref\Yb, \ref\Demaz) to check that
$$
\ga(X)Y_j\ga(X)^{-1}\ =\  X_j\de^{-(b_j,b_j)/2} Y_j \ =
\ \tau(Y_j).
$$
\proofbox

Actually, we can take here  an arbitrary
 $W$-invariant polynomial $g$ of $\{z_1,
\ldots, z_n\}$ such that
$b(\de^{g})\de^{-g}$ belong to $\C_\de[X]$ .

We claim that the
Schwartz functions $\{\ga p_b,\ b\in B_-\}$ defined
for the Macdonald
polynomials $\{p_b\}$ are pairwise orthogonal with respect to
the Fourier pairing $[\![\ ]\!]$. Here one should complete
\HH\ . Avoiding  this we will reformulate the statement
as follows:

\proclaim {Proposition}
The operators $\l_f^\ga\equal \ga \l_f\ga^{-1}$ defined for
 $f\in \C[y]^W$ are $W$-invariant
(see Theorem \ref\LF). Moreover, $\vph(\l_f^\ga)=\l_f^\ga$,
$\l_f^\ga (\ga p_b)=f(q^\rho\de^{-b})(\ga p_b)$, and the
corresponding eigenvalues (for all $f$)
distinguish different $\ga p_b$.
\endproclaim

%
%
%
%      REFERENCES
%
%
%
%\vskip 15pt
\AuthorRefNames [BGG]
\references
%\medskip
%\ninerm
%\baselineskip=11pt %!

[B]
\name{N. Bourbaki},
{\it Groupes et alg\`ebres de Lie}, Ch. {\bf 4--6},
Hermann, Paris (1969).

[C1]
\name{I.Cherednik},
{  Double affine Hecke algebras,
Knizhnik- Za\-mo\-lod\-chi\-kov equa\-tions, and Mac\-do\-nald's
ope\-ra\-tors},
IMRN (Duke M.J.) {  9} (1992), 171--180.

[C2]
\bibline,{ Double affine Hecke algebras and  Macdonald's
conjectures},
Annals of Mathematics {141} (1995), 191-216.

[C3]
\bibline,
{ Induced representations of  double affine Hecke algebras and
applications},
Math. Research Letters { 1} (1994), 319--337.

[C4]
\bibline,{ Difference-elliptic operators and root systems}, IMRN (1995).

[C5]
\bibline,
{ Integration of quantum many- body problems by affine
Knizhnik--Za\-mo\-lod\-chi\-kov equations},
Pre\-print RIMS--{  776} (1991),
(Advances in Math. (1994)).

[D]
\name{ C.F. Dunkl},
{ Hankel transforms associated to finite reflection groups},
Contemp. Math. {138} (1992), 123--138.

[EK1]
\name {P.I. Etingof}, and \name {A.A. Kirillov, Jr.},
{Macdonald's polynomials and representations of quantum
groups}, Preprint hep-th 9312103 (1993).

[EK2]
\bibline
{Representation-theoretic proof of the inner product and
symmetry identities for Macdonald's polynomials},
Compositio Mathematica (1995).

[GH]
\name{ A.M. Garsia}, and \name{M. Haiman},
{A graded representation model for Macdonald's polynomials},
Proc. Nat. Acad. Sci. USA {90}, 3607--3610.

[J]
\name{ M.F.E. de Jeu},
{The Dunkl transform }, Invent. Math. {113} (1993), 147--162.

[He]
\name{G.J. Heckman},
{  An elementary approach to the hypergeometric shift operators of
Opdam}, Invent.Math. {  103} (1991), 341--350.
omp. Math. {  64} (1987), 329--352.
v.Math.{  70} (1988), 156--236.

[KL]
\name{D. Kazhdan}, and \name{ G. Lusztig},
{  Proof of the Deligne-Langlands conjecture for Hecke algebras},
Invent.Math. {  87}(1987), 153--215.

[KK]
\name{B. Kostant}, and \name{ S. Kumar},
{  T-Equivariant K-theory of generalized flag varieties,}
J. Diff. Geometry{  32}(1990), 549--603.

[M1]
\name{I.G. Macdonald}, {  A new class of symmetric functions },
Publ.I.R.M.A., Strasbourg, Actes 20-e Seminaire Lotharingen,
(1988), 131--171 .

[M2]
\bibline, {  Orthogonal polynomials associated with root
systems},Preprint(1988).

[M3]
\bibline, {  Some conjectures for root systems},
SIAM J.Math. Anal.{  13}:6 (1982), 988--1007.

[O1]
\name{E.M. Opdam},
{  Some applications of hypergeometric shift
operators}, Invent.Math.{  98} (1989), 1--18.

[02]
\bibline, {Harmonic analysis for certain representations of
graded Hecke algebras},
Preprint Math. Inst. Univ. Leiden W93-18 (1993).

\endreferences

\bye